\documentclass[final]{siamltex}

\usepackage[dvips]{graphics}
\usepackage{epsfig}
\usepackage[left=3.3cm,top=2cm,right=3.2cm]{geometry}
\usepackage{amssymb}
\usepackage{amsmath}

\newtheorem{algorithm}[theorem]{Algorithm}

\title{Theoretical and experimental analysis of a Randomized
algorithm for Sparse Fourier transform analysis \thanks{This
        work was supported by NSF DMS-0245566, NSF DMS-0219233, AT\&T Research and DIMACS.}}

\author{Jing Zou\thanks{ Program in Applied and Computational
Mathematics, Princeton University, Princeton, NJ 08544
         {\tt jzou@math.princeton.edu}.}
        \and Anna Gilbert \thanks{ Department of Mathematics, University of Michigan.}
        \and Martin Strauss \thanks{ Departments of Mathematics and Electrical Engineering and Computer Science, University of Michigan.}
\and Ingrid Daubechies \thanks{ Program in Applied and Computational
Mathematics and Department of Mathematics, Princeton University, Princeton, NJ 08544.}}

\begin{document}

\maketitle

\begin{abstract}
We analyze a sublinear RA$\ell$SFA (Randomized Algorithm for Sparse
Fourier Analysis) that finds a near-optimal $B$-term Sparse Representation $R$ for a
given discrete signal $S$ of length $N$, in time and space $poly(B,\log(N))$, following the approach given in \cite{GGIMS}.
Its time cost $poly(\log(N))$ should be compared with the superlinear
 $\Omega(N \log N)$ time requirement of the Fast
Fourier Transform (FFT). A straightforward implementation of the RA$\ell$SFA, as presented in the
theoretical paper \cite{GGIMS}, turns out to be very slow in practice. Our main result is a greatly
improved and practical RA$\ell$SFA. We introduce several new ideas and techniques that 
speed up the algorithm. Both rigorous and heuristic arguments for parameter choices are presented.
Our RA$\ell$SFA constructs, with probability at least $1-\delta$, a near-optimal $B$-term representation
$R$ in time $poly(B)\log(N)\log(1/\delta)/\epsilon^{2} \log(M)$ such that $\Vert S-R\Vert^{2}_2\leq(1+\epsilon)\Vert
S-R_{opt}\Vert^{2}_2$. Furthermore, this RA$\ell$SFA implementation already beats the FFTW for not unreasonably large $N$. We
extend the algorithm to higher dimensional cases both theoretically and numerically. The crossover point lies at
$N \simeq 70000$ in one dimension, and at $N \simeq 900$ for data on a $N \times N$ grid in two dimensions for
small $B$ signals where there is noise.
\end{abstract}



\pagestyle{myheadings}
\thispagestyle{plain}
%\markboth{June 2}{An Improved Randomized algorithm for Sparse Fourier transform analysis}
%\markboth{Sept. 8}{Sept. 8}
\markboth{J. ZOU, A. GILBERT, M. STRAUSS, I. DAUBECHIES}{Theoretical and Experimental Analysis of a Randomized
Algorithm for Sparse Fourier Analysis}

\section{Introduction}

We shall be concerned with discrete signals $S=(S(0),\ldots, S(N-1)) \in \mathbb{C}^N$ and their Fourier
transforms $\hat{S}=(\hat{S}(0),\ldots,\hat{S}(N-1))$, defined by
$\hat{S}(\omega)=\frac{1}{\sqrt{N}}\sum_{t=0} ^{N-1} S(t)e^{-2 \pi i \omega t/N}$. In terms of the Fourier basis
functions $\phi_{\omega}(t) = \frac{1}{\sqrt{N}} e^{2 \pi i \omega t/N}$, $S$ can be written as $S
= \sum_{\omega=0}^{N-1} \hat{S}(\omega)\phi_{\omega}(t)$; this is the (discrete) Fourier representation of $S$.

In many situations, a few large Fourier coefficients already capture the major time-invariant wave-like information
of the signal and very small Fourier coefficients can thus be discarded. The problem of finding the (hopefully few)
largest Fourier coefficients of a signal that describe most of the signal trends, is a fundamental task in
Fourier Analysis. Techniques to solve this problem are very useful in data compression, feature extraction, finding
approximating periods and other data mining tasks \cite{GGIMS}, as well as in situations where multiple scales
exist in the domain (as in e.g. materials science), and the solutions have sparse modes in the frequency domain.

Let $S$ be a signal that is known to have a sparse $B$-term Fourier
representation with $B\ll N$, i.e.,
\begin{equation}\label{Brepn}
S(t)=\frac{1}{\sqrt{N}}(a_{1}e^{i2\pi\omega_{1}t/N}+\ldots+a_{B}e^{i2\pi\omega_{B}t/N}),
\end{equation}
and let us assume that it is possible to evaluate $S$, at arbitrary $t$, at cost $O(1)$ for every evaluation.

To identify the parameters $a_{1},\ldots,a_{B},\omega_{1},\ldots,\omega_{B}$, one can use the Fast Fourier
Transform (FFT). Starting from the $N$ point-evaluations $S(0), \ldots, S(N-1)$, the FFT computes all the
Fourier coefficients; one can then take the largest $B$ coefficients and the corresponding modes. The time cost
for this procedure is $\Omega (N \log N)$; this can become very expensive if $N$ is huge. (Note that
all logarithms in this paper are with base 2, unless stated otherwise.) The problem becomes worse in higher dimensions. If one uses grids of size $N$ in
each of $d$ dimensions, the total number of points is $N^{d}$ and the FFT procedure takes $\Omega(d N^d \log N)$ time.
It follows that
 identifying a sparse number of modes and amplitudes is expensive for even fairly modest $N$.
Our goal in this paper is to discuss much faster algorithms that can identify the coefficients $a_{1},\ldots,a_{B}$ and
the modes $\omega_{1},\ldots,\omega_{B}$ in equation (\ref{Brepn}). These algorithms will not use all the samples
$S(0),\ldots, S(N-1)$, but only a very sparse subset of them.

In fact, we need not restrict ourselves to signals that are exactly equal to a $B$-term representation. Let us
denote the optimal $B$-term Fourier representation of a signal $S$ by $R_{opt}^B(S)$; it is simply a truncated version of the Fourier representation of $S$, retaining only the $B$ largest coefficients.
We are then interested in identifying (or finding a close approximation to) $R_{opt}^B(S)$ via a fast algorithm.
The papers \cite{GGIMS} \cite{Mansour} \cite{GMS} provide such algorithms; all compute a (near-)optimal $B$- term Fourier
representation $R$ in time and space $poly(B,\log(1/\delta),\ \log N,\log M,1/\epsilon)$, such that $\Vert
S-R\Vert_2^{2}\leq(1+\epsilon)\Vert S-R_{opt}^B(S)\Vert_2^{2}$, with success probability at least $1-\delta$,
where $M$ is an a priori given upper bound on $\|S\|_2$. The
algorithms in these papers share the property that they need only some random subsets of the input rather than all the data;
they differ in many details: the different papers assume different conditions on $N$, (for example, $N$ is assumed to be
a power of 2 or a small prime number in \cite{Mansour}; N may be arbitrary
but is preferably a prime in \cite{GGIMS}); the algorithms also use different
schemes to locate the significant modes. (Here we say a mode $\omega$ is significant if for
some pre-set $\eta$, $|\hat{S}(\omega)|^2 \geq \eta \|S\|^2$.) Mansour and Sahar \cite{Sahar} implemented a similar algorithm
for Fourier analysis on the set $\mathbb{Z}_2^n$, where our algorithm is for Fourier analysis on $\mathbb{Z}_N$.

The results of \cite{GGIMS} can be extended to more general representations, with respect to a particular basis or a family of bases; examples are wavelet
bases, wavelet packets or Fourier bases. We shall use the acronym RA$\ell$STA (Randomized Algorithm for Sparse
Transform Analysis) for this family of algorithms. We here restrict ourselves to the Fourier case and thus RA$\ell$SFA.

For a wide range of applications, the speed potential suggested by the sublinear cost of these
algorithms is of great importance. In this paper, we concentrate on the approach proposed in \cite{GGIMS}.
Note that \cite{GGIMS} gives a theoretical rather than a practical analysis in the sense that it does not discuss parameter settings;
it gives few hints about the order of the polynomial in $B$ and $\log N$; in fact, a straightforward implementation of RA$\ell$SFA following the set-up of \cite{GGIMS} turns out to be too
slow to be practical, so that none of the direct implementation work was published. In addition, \cite{GGIMS} did not discuss extensions
 to higher dimensions, where the pay-off of RA$\ell$SFA versus the FFT is expected to be larger.

Our main result in this paper is a version of RA$\ell$SFA that addresses these problems. We give
theoretical and heuristic arguments for the setting of parameters; we introduce some new ideas that produce a practical RA$\ell$SFA implementation. Our new version can outperform the FFTW when $N$ is around $70,000$ and $B$ is small.

\textbf{ A Motivating Example.} RA$\ell$SFA is an exciting replacement for the FFT to solve multiscale models.
Typically, one wants to simulate a multiscale model in several dimensions with both a microscopic and a
macroscopic description. The solution to the model has rapidly oscillating coefficients with period proportional to a small parameter
$\epsilon$. For examples of multiscale problems of size $N$ that are dominated by the behavior of $B\ll N$ Fourier components, see e.g \cite{BLP}. In a traditional (pseudo-)spectral method, one computes the spatial derivatives by the FFT
and Inverse FFT at each time iteration; consequently the time to find the Fourier representation of a signal is
the determining factor in the overall time of simulation. In multiscale problems,
where only a small number of Fourier modes contribute to the energy of an initial condition and coefficient
functions, we expect that RA$\ell$SFA will
significantly speed up the calculation for large $N$. In fact, a preliminary study has shown \cite{OLOF}
that for some transport and diffusion equations with multiple scales, using only significant frequencies to approximate intermediate solutions does not substantially degrade the quality of the
approximate final solution to the multiscale problem. By using the most significant frequencies and RA$\ell$SFA instead
of all frequencies and the FFT, we could replace a superlinear algorithm by a poly-log (polynomial in the logarithm)
algorithm.  The corresponding decrease of the running time would make it possible to handle a larger number of grid points in
high dimensions. We shall present detailed applications
of this algorithm in multiscale problems in \cite{ZDR}.

\textbf{Notation and Terminology.} For any two frequencies $\omega_1$, $\omega_2$,
where $\omega_1 \neq \omega_2$, we say that $\hat{S}(\omega_1)$ is bigger than
$\hat{S}(\omega_2)$ if $| \hat{S}(\omega_1)|>|\hat{S}(\omega_2)|$. The squared norm $\|S\|_2^2=\sum_{t=0}^{N-1}|S(t)|^2$
of $S$ is also called the energy of $S$; we shall refer to
$| \hat{S}(\omega) |^2$ as the energy of the Fourier coefficient $\hat{S}(\omega)$. Similarly, the energy of a
set of Fourier coefficients is the sum of the squares of their magnitudes. We shall use only the $\ell ^2$-norm in this
paper; for convenience, we therefore drop the subscript from now on, and denote $\|F\|_2^2$ by $\|F\|^2$
for any signal $F$.

We denote the convolution by $F*G$, $(F*G)(t)=\sum_s{F(s)G(t-s)}$. It follows that $\widehat{F*G}=\sqrt N \hat{F}\hat{G}$.
We denote by $\chi_T$ the
signal that equals 1 on a set $T$ and zero elsewhere. The index to $\chi_T$ may be either time or frequency; this is
made clear from context. For more background on Fourier analysis, see \cite{weaver}. The support $supp(F)$ of
a vector $F$ is the set of $t$ for which $F(t)\neq 0$. A signal is 98$\%$ pure if there exists a frequency $\omega$ and some signal $\rho$, such
that $S=a \phi_{\omega}+\rho$ and $|a|^2\geq 0.98\|S\|^2$.

RA$\ell$SFA is a randomized algorithm. By this, we do {\bf not} mean the signal is randomly chosen from some kind of
distribution, with our timing and memory requirement estimates holding
with respect to this distribution; on the contrary, the signal, once given to us, is {\bf fixed}.
 The randomness lies in the algorithm. After random sampling, certain operations are
repeated many times, on different subsets of samples, and averages and medians of the results are computed.
 We set in advance
a desired probability of success $1-\delta$, where $\delta>0$ can be arbitrarily small. Then the claim is that
for each arbitrary input $S$, the algorithm succeeds with probability $1-\delta$,
 i.e., gives a $B$-term estimate $R$ such that
$\|S-R\|^2 \leq (1+\epsilon) \|S-R_{opt}^B\|^2$. For given $\epsilon$, $\delta$, numerical experiments show that the
algorithm may take $O(B^2 \log N)$ time and space.

\textbf{Organization.} The chapters are organized as follows. Section 2 shows the testbed and numerical experiments about the
comparison of our RA$\ell$SFA and the FFTW. In Section 3, we
introduce all the new techniques and ideas of RA$\ell$SFA (different from \cite{GGIMS}) and its extension to multi-dimensions.

\section{Testbed and Numerical Results of RA$\ell$SFA}

In this section, we present numerical results of RA$\ell$SFA. We begin in Section
\ref{subset:onedim} with comparing the running time of RA$\ell$SFA and the FFTW
for some one dimensional test examples. In Section \ref{sect:ntwodim},
the performances of two dimensional RA$\ell$SFA and the FFTW for some test signals
are shown.

The randomness of the algorithm implies that the performance differs each time for
the same group of parameters. Hence, we give the average data, bar and quartile
graph based on 100 runs as well as the fastest data among these experiments. The
popular software FFTW \cite{FJ} version 2.1.5 is used to determine the
timing of the Fast Fourier Transform for the same data.

The test signals are either superpositions of $B\ll N$ modes in the frequency domain, that is, $S=\sum_{j=1}^{B} c_j
\phi_{\omega_j}$, contaminated with Gaussian white noise, or signals for which the Fourier coefficients exhibit rapid decay, so that a $B$-mode approximation with $B\ll N$ will already be very accurate. Different choices of the $\omega_j$ were checked; these did not influence the whole execution time. These choices included cases where some frequencies were close; note that this is the ``hard'' case for most
estimation algorithms. For RA$\ell$SFA, which contains random scrambling
operations (that are later described), the distance between the modes does not
matter if $N$ is prime. If $N$ is not prime, then $gcd(\omega_1-\omega_2, N)$
cannot decrease by the scrambling operation, so that different $(\omega_1,
\omega_2)$ pairs may (in theory) lead to different performances; in practice,
this doesn't seem to matter. In all these situations,
RA$\ell$SFA reliably estimates the size and locations of the few largest coefficients. We also set other parameters as follows: accuracy factor
$\epsilon=10^{-2}\|S\|$, failure probability $\delta=0.05$.

The parameter choices in the algorithm are quite tricky. The
theoretical bounds given in \cite{GGIMS} do not work well in practice; instead
much smaller parameters and heuristic settings work more efficiently.

All the experiments were run on an AMD Athlon(TM) XP1900+ machine with Cache size 256KB, total memory 512 MB, Linux kernel version 2.4.20-20.9 and compiler gcc version 3.2.2.

\subsection{Numerical Results in one dimension}
\label{subset:onedim}

The first implementation results of RA$\ell$SFA were not published; the program was
basically a proof of concept, not optimized. With the choices and parameters
described in \cite{GGIMS}, it was extremely slow and thus not practical for
real-world applications. The implementation we present here runs several order
of magnitude faster; this involves introducing many adjustments and ideas to
the algorithm of \cite{GGIMS}. (See Section 3 for details.)

The goal of this paper is to check the possibility to replace the FFT with RA$\ell$SFA for sparse and long signals. Therefore, we focus on comparing the performance of RA$\ell$SFA and FFTW in the following subsections.

\subsubsection{Experiments for an Eight-mode Representation}
\label{subset:eightoned}

We begin with the experiments for recovering a signal consisting of eight modes (with and without noise).
In the noisy signal case, the noise is a Gaussian white noise with signal-to-noise ratio
($SNR$, defined as $10\log_{10}\frac{\|S\|^2}{N\sigma^2}$) approximately 5dB.  The coefficients are randomly taken
from the interval $[1,10]$ and the significant modes from $[0, N-1]$.

Two kinds of running time for each algorithm are provided. One is the total running time and
another is the running time excluding the sampling time. As we know, the FFT
takes $\Omega(N)$ to compute all signal values. On the other hand, our algorithm
doesn't need all the sample values. All our conclusions are based on the time
{\it excluding} the sampling. However, we still list the running time including sampling time
as well because of the existence of various forms of data in practice. For example, in pseudospectral
applications, the data need to be computed from a
B-superposition, which may take $O(B)$ per sample. It is possible to sample more quickly, which
 is addressed in \cite{GMS}. On the other hand, if the data is
already stored in a file or a disk, we simply get them without any computation.
In all these cases, we assume the data is either already in memory or
available through computation. Thus we don't need to go through every data, which would take time $O(N)$.

Table \ref{tab:B8onedim} provides a
comparison of the running times of the FFTW and RA$\ell$SFA for eight-mode clean and noisy signals.
 In the beginning when $N$ is small, the FFTW is almost instantaneous. As
the signal length $N$ increases, its time grows superlinearly. On the contrary, RA$\ell$SFA takes longer
 time in smaller
$N$ cases; however the time cost remains almost constant regardless of the
signal length. In addition, the benchmark FFTW
software fails to process more than $10^8$ data because it runs out of the memory space. In contrast,
RA$\ell$SFA has no difficulty at all since it does not need all the data. A simple interpolation from the entries in Table 
\ref{tab:B8onedim} predicts that RA$\ell$SFA beats the FFTW when $N>15,200$ for eight-mode signals, all the more convincingly when $N$ is larger.
  If we compare the time including sample computation, the cross-over point would be $N=70,000$. The table also
shows the linear relationship between the time cost and the logarithm of the length $N$.

\begin{table}[htbp]
\begin{center}
\begin{tabular}{|l|l|l|l|l|l|l|}
\hline
Length& \multicolumn{2}{l|}{Time of}&  \multicolumn{1}{l|}{ Time of} & \multicolumn{2}{l|}{Time of RA$\ell$SFA}&\multicolumn{1}{l|}{Time of FFTW}\\
N &  \multicolumn{2}{l|}{RA$\ell$SFA}  &  \multicolumn{1}{l|}{FFTW} &\multicolumn{2}{l|}{(excluding sampling) }&
\multicolumn{1}{l|}{(excluding sampling)}\\ \cline{2-6}
    &clean& noisy&  &clean& noisy& \\ \hline \hline
$10^{3}$&
 0.22 & 0.25 & 0  &0.01 &0.02 &
0 \\ \hline $10^{4}$&
 0.25&
 0.29 & 0.04& 0.03 & 0.04&
0.01
 \\ \hline $10^{5}$&
 0.32 &
 0.34&  0.46&0.05  & 0.05&
 0.17 \\ \hline
$10^{6}$&
 0.37 & 0.41&5.01 & 0.07&
 0.08 &
2.23\\ \hline $10^{7}$&
 0.44 & 0.48& 54.57  & 0.10&
 0.11 &
26.24\\ \hline
\end{tabular}
\end{center}
\caption{Time Comparison between RA$\ell$SFA and FFTW (B=8) based on 100 runs. ``Clean'' means that the test signal is pure. ``Noisy'' means the signal is contaminated with noise of $SNR=5dB$. ``Excluding Sampling'' column lists the running time without precomputation of sample values.} \label{tab:B8onedim}
\end{table}   \addtocounter{figure}{+1}

\begin{figure}[htbp]
\begin{center}
\includegraphics[%
  width=6cm]{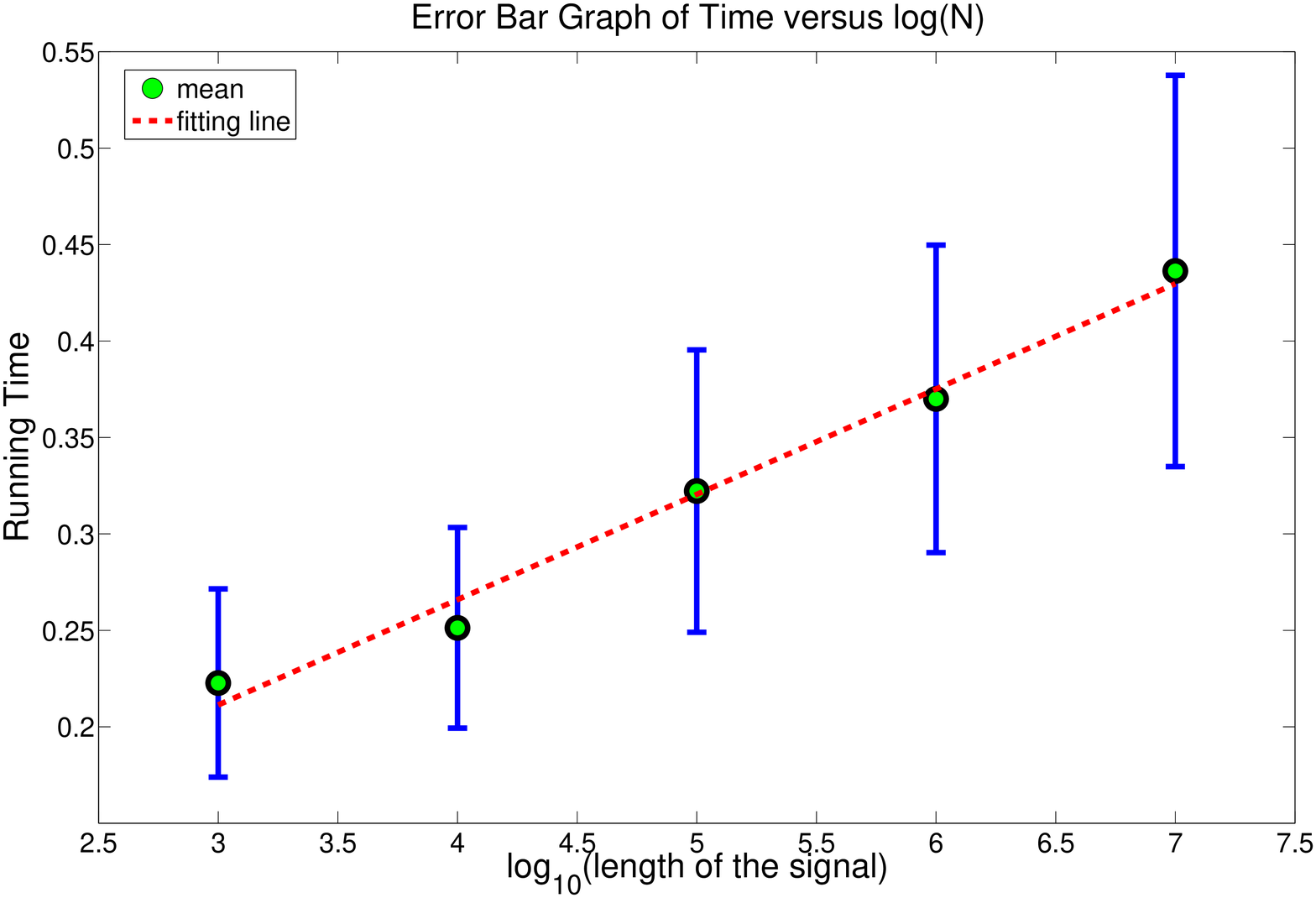}
\includegraphics[%
  width=6cm]{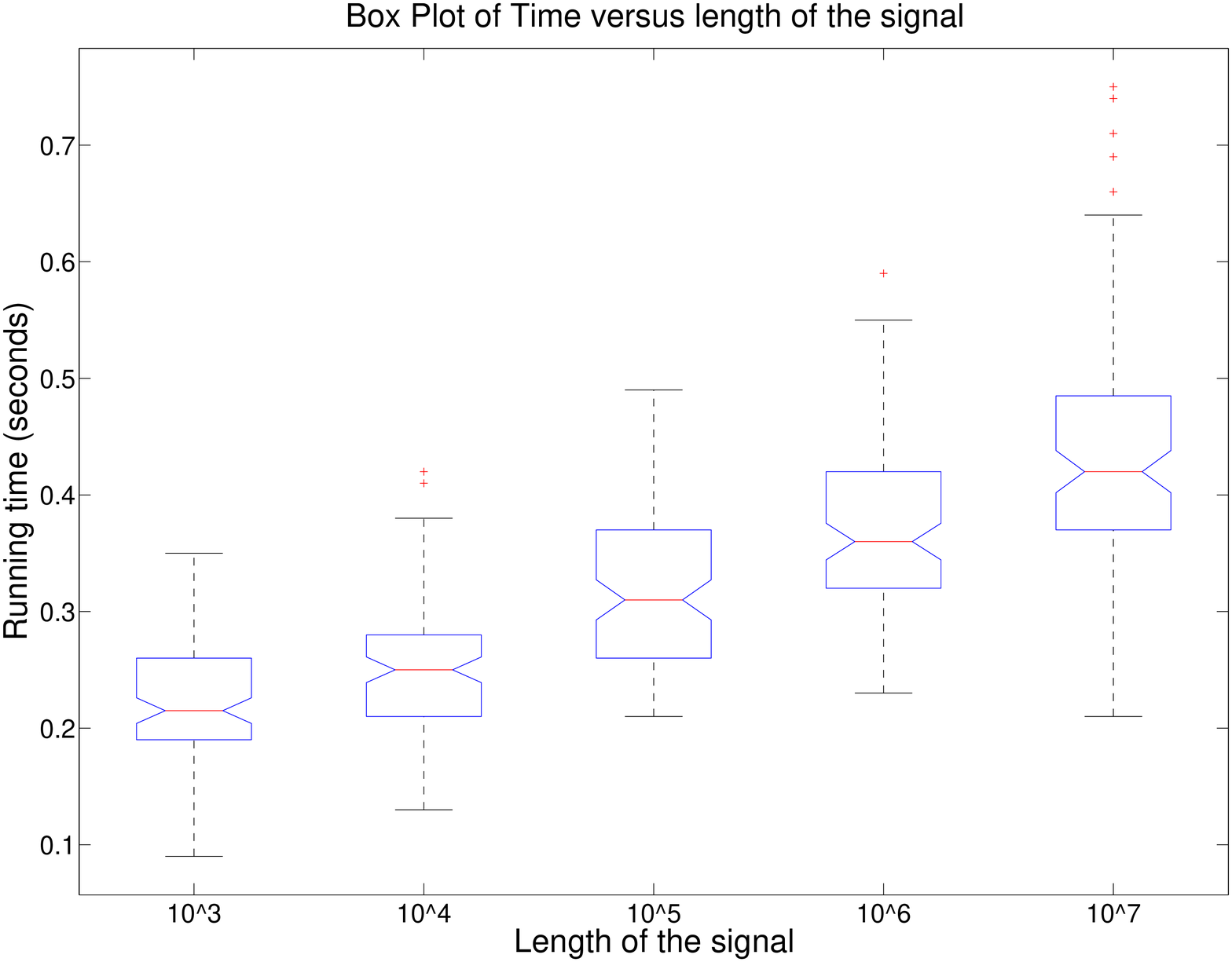}
\caption{Comparing the total running times of 8-mode RA$\ell$SFA for 100 different runs of the randomized algorithm. Left: mean and
variance as a function of $N$; right: median, quartiles and total spread of the runs as a function of $N$, $B=8$}\label{fig:B8bar1}
\end{center}
\end{figure}

As can be expected from a randomized algorithm, RA$\ell$SFA has a different performance in each run. Figure \ref{fig:B8bar1}
illustrates the spread of the execution time (including sampling) for pure signals over 100 runs. 

    \addtocounter{table}{+1}

\subsubsection{Experiments with Different Levels of Noise}
\label{subset:noise}
\begin{figure}[htbp]
\begin{center}
\includegraphics[%
  width=6cm]{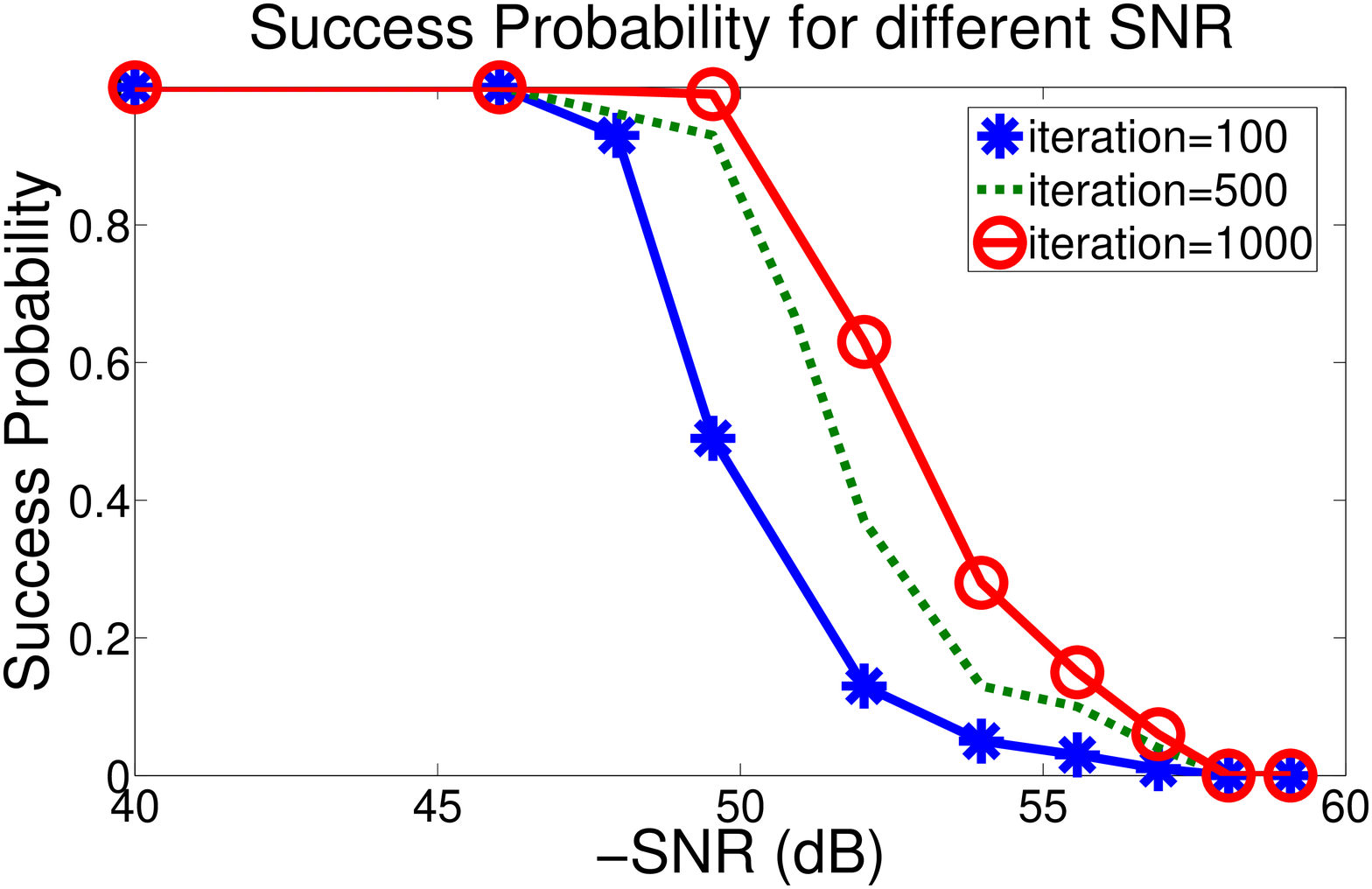}
\includegraphics[%
  width=6cm]{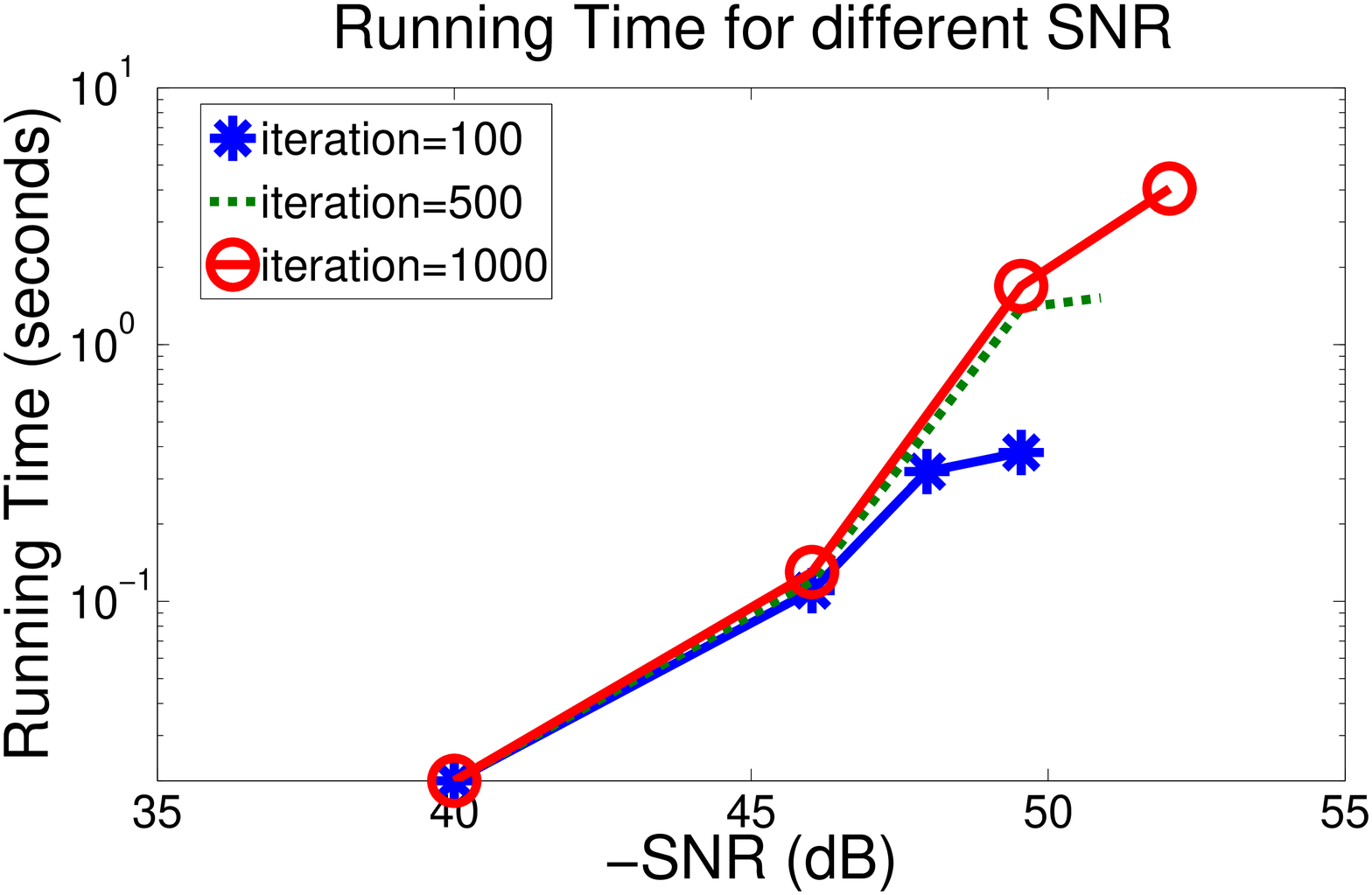}
\caption{Experiments for signal $S=\phi_0+noise$ with length $N=10,009$. Compare the success rate 
and running time of RA$\ell$SFA when the total number of iterations is bounded by 100, 500 and 1000 (respectively), based on 100 different runs of the randomized algorithm
in each case. Left: success probability of RA$\ell$SFA as $SNR$ decreases; right: running time of RA$\ell$SFA as $SNR$ decreases (we only show the running time when success probability is greater than $50\%$). Note that the abscissa show $-SNR$ each time, meaning that the $\ell^2$-norm of the noise is much larger than that of the signal in the regimes
illustrated here; for instance, $SNR=-60dB$ means that the $\ell^2$-norm of the noise equals $1,000 \times \ell^2 \,norm \,of \,the \,signal $.} \label{fig:noise}
\end{center}
\end{figure}

In the experiments above, we compared the performance of clean and slightly noisy signals. Here, we shall push
the noise level much higher, keeping $N$ and $B$ fixed to illustrate the effect of noise. 
Also, instead of allowing the algorithm to run for $poly(B, \log N, 1/\epsilon, \log(1/\delta) )$ iterations,
we set a smaller fixed upper bound (so that the success probability is no longer $1-\delta$). 
When noise is present, 
it influences the success probability with which modes with small amplitude are detected. To explore this, we ran an experiment 
with only a single mode; we kept the amplitude of the mode constant and increased the noise. Figure \ref{fig:noise} (left) shows the \textit{success probability} of the detection of the single mode by the algorithm (estimated
by running 100 trials each time and recording the number that were successful) for three different settings of the maximum number of iterations. 

The dependence of the \textit{running time} on the $SNR$ in the case of detection of a single mode is illustrated in
Figure \ref{fig:noise} (right), where we show the results of the average over 100 runs for every data point, with only a very loose a priori 
restriction on the running time ($\leq$1000 iterations); only parameter settings with over $50\%$ success probability were
taken into account.

    \addtocounter{table}{+1}

This experiment indicates
that it is possible to detect modes that are significantly weaker than the noise, within limits, of course. If the amplitude
of the signal is too weak, then trying to detect it may waste many resources. In practice we
shall put our cut-off on the amplitude at about one sixth of the noise level, i.e., at $\sigma/6$; this can of course be adjusted
depending on whether one wishes fast speed or not.

Although $SNR$ is the standard characterization of noise intensity, it is not clear that it is the parameter that 
matters most for our algorithm. We therefore also ran an experiment in which we compare the results for two 
different values of $N$: 10,009 (as in the Figures above) and 100,003, respectively. 
The second value of $N$ is about 10 times larger than the first; for the same choices of $\sigma$ and $c$ (the amplitude of the single mode), the $SNR$ 
for the second $N$ is smaller by $10dB$. Table \ref{tab:noisesnr}, comparing the performance for these two values of $N$ and several choices of 
$\sigma$, shows that the value of $\sigma$ itself rather than $SNR$ governs the running time and success probability.

\begin{table}[htbp]
\begin{center}
\begin{tabular}{|c||c|c|c||c|c|c|}
\hline
   & \multicolumn{3}{c||}{$N_1=10,009$}&  \multicolumn{3}{c|}{$N_2=100,003$}\\ \hline 
$\sigma$& success probability & time & $SNR$ & success probability & time & $SNR$ 
 \\ \hline 2& $100\%$ & 0.11& -46.02& $100\%$&0.19 & -56.02\\ \hline
 2.5 & $93\%$ & 0.32& -47.96& $77\%$& 0.55& -57.96\\ \hline 3& $49\%$ & 0.38& -49.54& $27\%$&0.61 &-59.54 \\ \hline 3.5& $21\%$ & 0.45& -50.88& $10\%$& 0.38&-60.88 
\\ \hline 4& $13\%$ & 0.38& -52.04& $1\%$& 0.37&-62.04 \\ \hline
\end{tabular}
\end{center}
\caption{Exploring the dependence on $\sigma$ versus $SNR$ of the influence of the noise on processing
the signal $S=\phi_0+noise$, where the noise is gaussian $N(0, \sigma)$. For two different values of $N$,
$N_1=10,009$ and $N_2=100,003 \approx 10N_1$, respectively, and a range of values for $\sigma$, we determined the success
probability within 100 runs, and the average running time for successful runs. In both cases we see a clear transition as 
$\sigma$ increases; the location of the transition (between 2.5 and 4 for $N_1$, between 2 and 3.5 for $N_2$) shifts slightly 
with $N$, but it is nevertheless clear that $\sigma$ is a better parameter to track than $SNR$: in fact, the largest choice for 
$\sigma$, $\sigma=4$, still has lower $SNR$ in the case $N=N_1$ than the smallest choice, $\sigma=2$, for $N=N_2$, yet the success
probability and running time are much worse. } \label{tab:noisesnr}
\end{table}   \addtocounter{figure}{+1}

\subsubsection{Experiments with Different Numbers of Modes}
\label{sect:diffB}

The crossover points for $N$ are different for signals with different $B$; the number of modes has an important
influence on the running time. To investigate this, we experimented with fixed $N$ (we took a prime number
$N=2,097,169$ (a prime number) for RA$\ell$SFA and $N=2^{21}=2,097,152$ for FFTW) but varying $B$. In all cases, we take $S$ to be
a superposition of exactly $B$ modes, i.e., $S(t)=\sum_{i=1}^{B}{c_{i}\phi_{\omega_i}}$ for some $B$. Table
\ref{tab:diffB} compares the running time for different $B$ using the FFTW and RA$\ell$SFA. For small $B$, RA$\ell$SFA
takes less time because $N$ is so large. The execution time for the FFT can be taken to include the time for
evaluation of all the samples (which increases linearly in $B$) or not (in which case the execution time is
constant to $B$). In both cases, the FFTW overtakes RA$\ell$SFA as $B$ increases; the execution time of the FFTW
is constant or linear in the number of modes $B$ (depending on whether the evaluation of samples is included),
while that of RA$\ell$SFA is polynomial of higher order. For $N=2,097,169$, the FFTW is faster than RA$\ell$SFA
when $B \geq 33$. By regression techniques on the experimental data, one empirically finds that the order of $B$
in RA$\ell$SFA is quadratic. This is the main disadvantage of RA$\ell$SFA. (Although this
nonlinearity in $B$ was expected by the authors of \cite{GGIMS}, the observation that it played such an
important role even for modest $B$ was the motivation for Gilbert, Muthukrishnan and Strauss to construct in
\cite{GMS} a different version of RA$\ell$SFA that is linear in $B$ for all $N$.) Hence, RA$\ell$SFA is most
useful for a long signal with a small number of modes.

\begin{table}[htbp]
\begin{center}
\begin{tabular}{|c|c|c|c|c|}
\hline
Number of modes & Time of & Time of & Time of RA$\ell$SFA &  Time of FFTW\\
B & RA$\ell$SFA & FFTW & (exclude sampling )& (exclude sampling)\\ \hline
 $2$&
 0.05 &
  7.49&
 0.03 &
5.46\tabularnewline
$4$&
 0.14 &
 9.38&
 0.05 &
5.46\tabularnewline
$8$&
 0.35 &
 13.22 &
0.07 &
5.46\tabularnewline
$16$&
 2.48 &
 20.92 &
0.83 &
5.46\tabularnewline
$32$&
 15.53 &
 36.28 &
4.13 &
5.46\tabularnewline
$64$&
 107.55 &
 67.16 &
39.55&
5.46  \tabularnewline
\hline
\end{tabular}
\end{center}
\caption{Time Comparison between RA$\ell$SFA and FFTW for Different $B$ when $N \approx 2,097,169$} \label{tab:diffB}
\end{table}   \addtocounter{figure}{+1}

\subsubsection{Experiments with Signals that have infinitely many modes with rapid decay in frequency}

For our final batch of one-dimensional experiments, we ran the algorithm on the signal $S=1/(1.5+\cos 2\pi t)+noise$. In continuous time, the clean signal has infinitely many modes with amplitudes that decay
exponentially as the frequency of the mode increases. We ran the experiment with a white Gaussian noise once with $SNR$
$-20dB$ and a second time with $SNR=-8dB$, with $N=1000$. The threshold for the amplitudes of modes we wished to find was adjusted to the
noise level in both cases.

\begin{figure}[htbp]
\begin{center}
\includegraphics[%
  width=7cm]{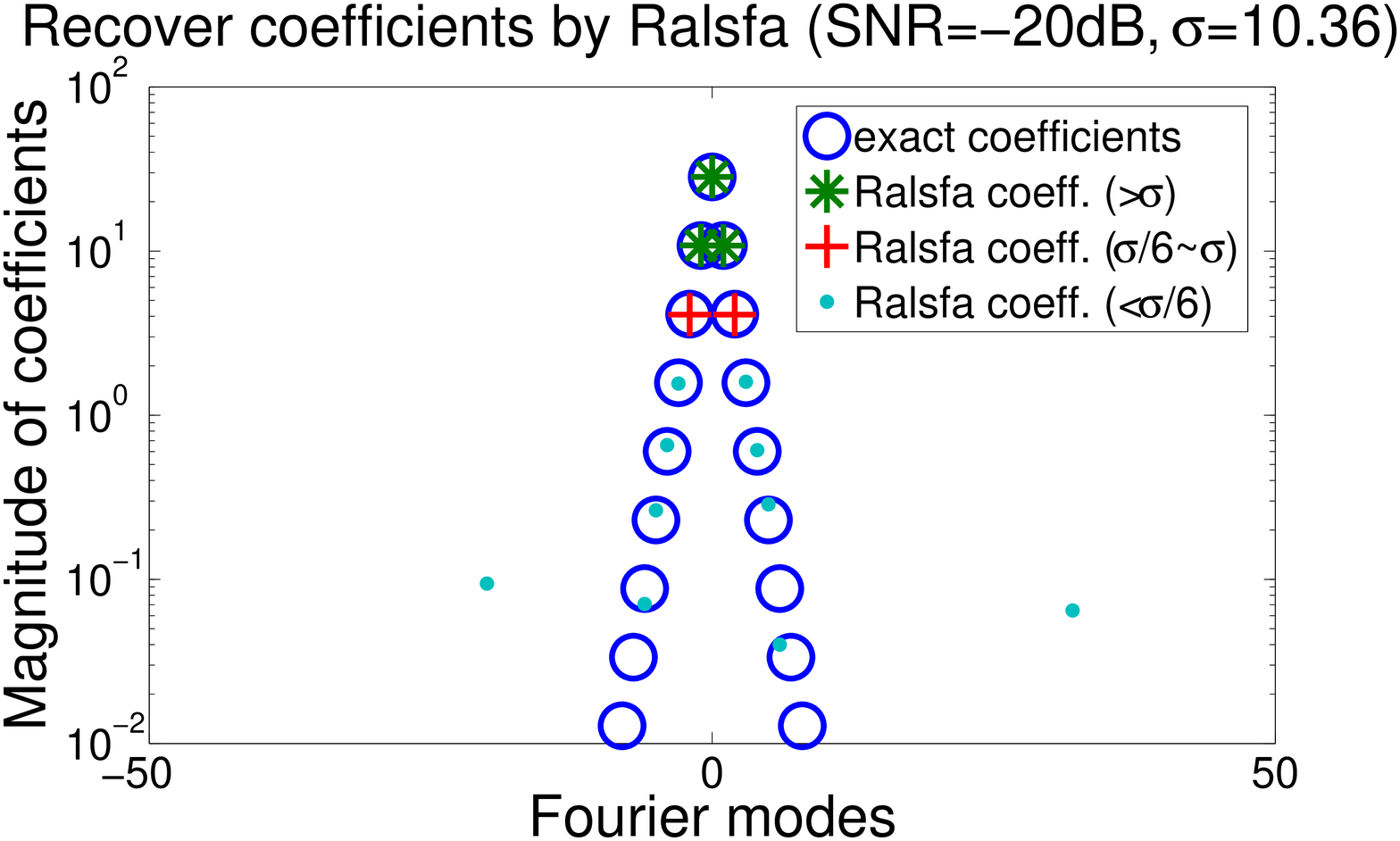}
\includegraphics[%
  width=7cm]{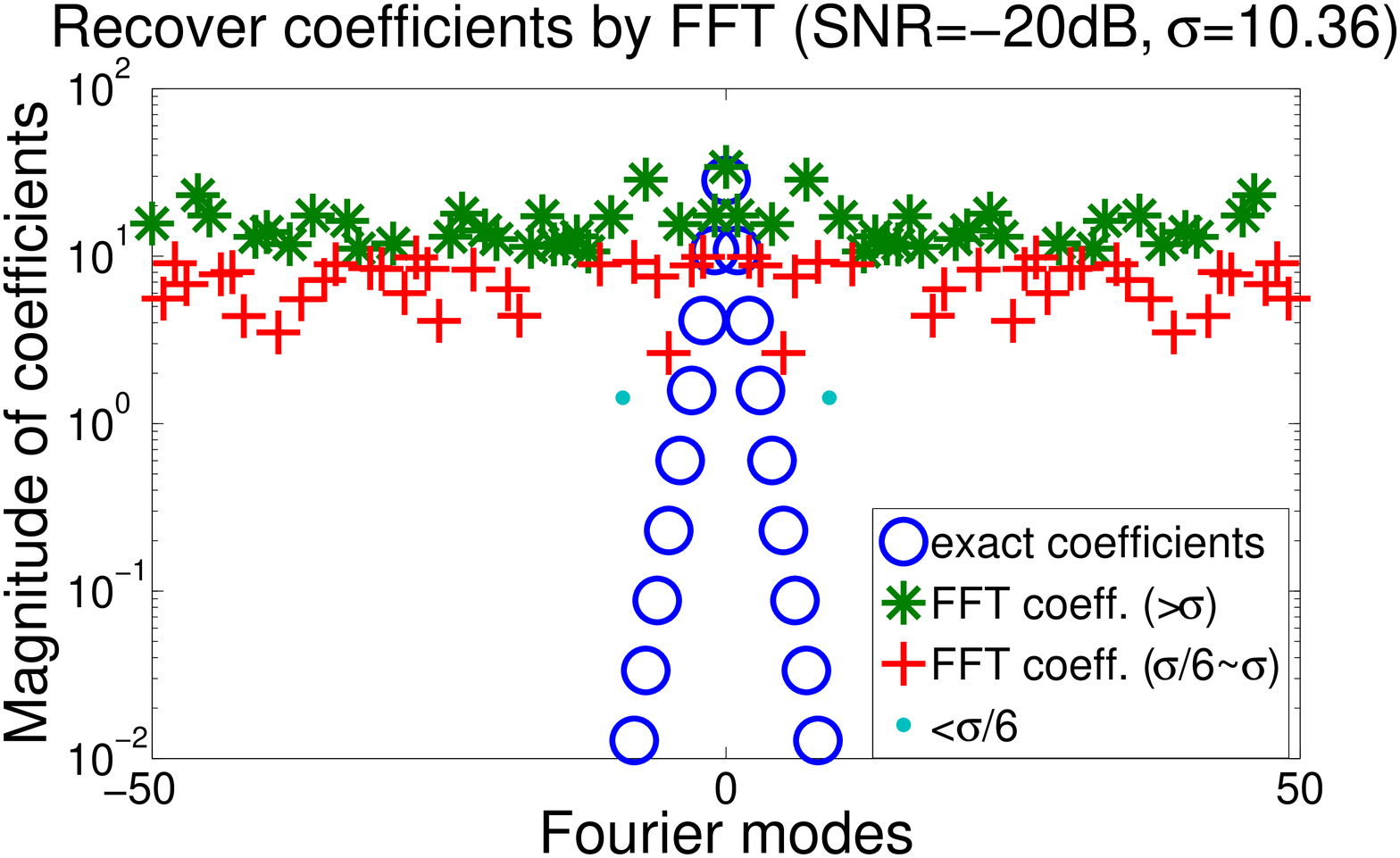}
\caption{For signal $S=1/(1.5+\cos 2\pi t)+noise$ with $SNR=-20dB$. Compare the approximation effect by RA$\ell$SFA and FFTW. Left: 
approximation of the significant coefficients 
by RA$\ell$SFA; the relative approximation error is $0.74\%$; right: approximation of the significant coefficients 
by FFTW. } \label{fig:decay}
\end{center}
\end{figure}    \addtocounter{table}{+1}

\begin{figure}[htbp]
\begin{center}
\includegraphics[%
  width=7cm]{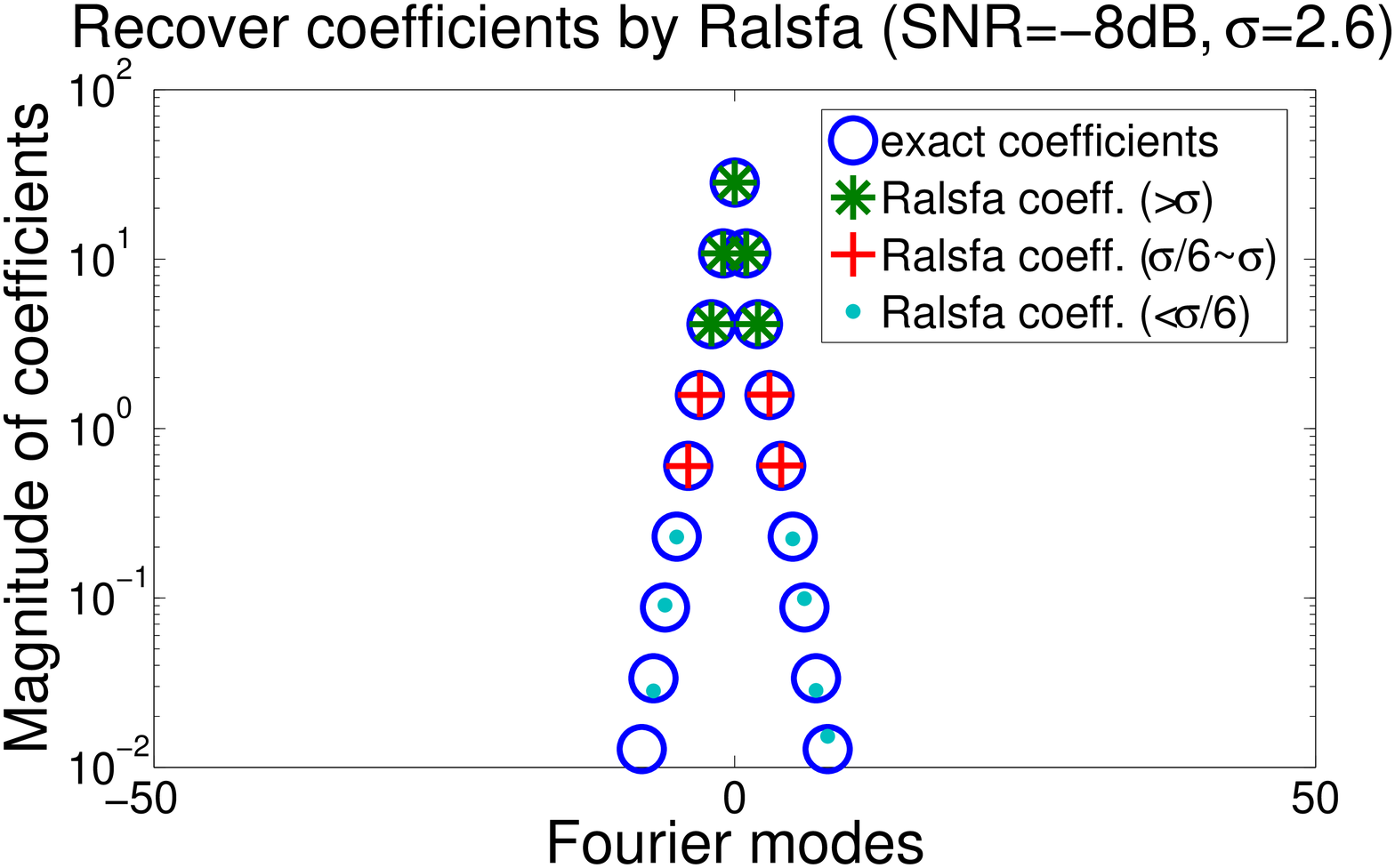}
\includegraphics[%
  width=7cm]{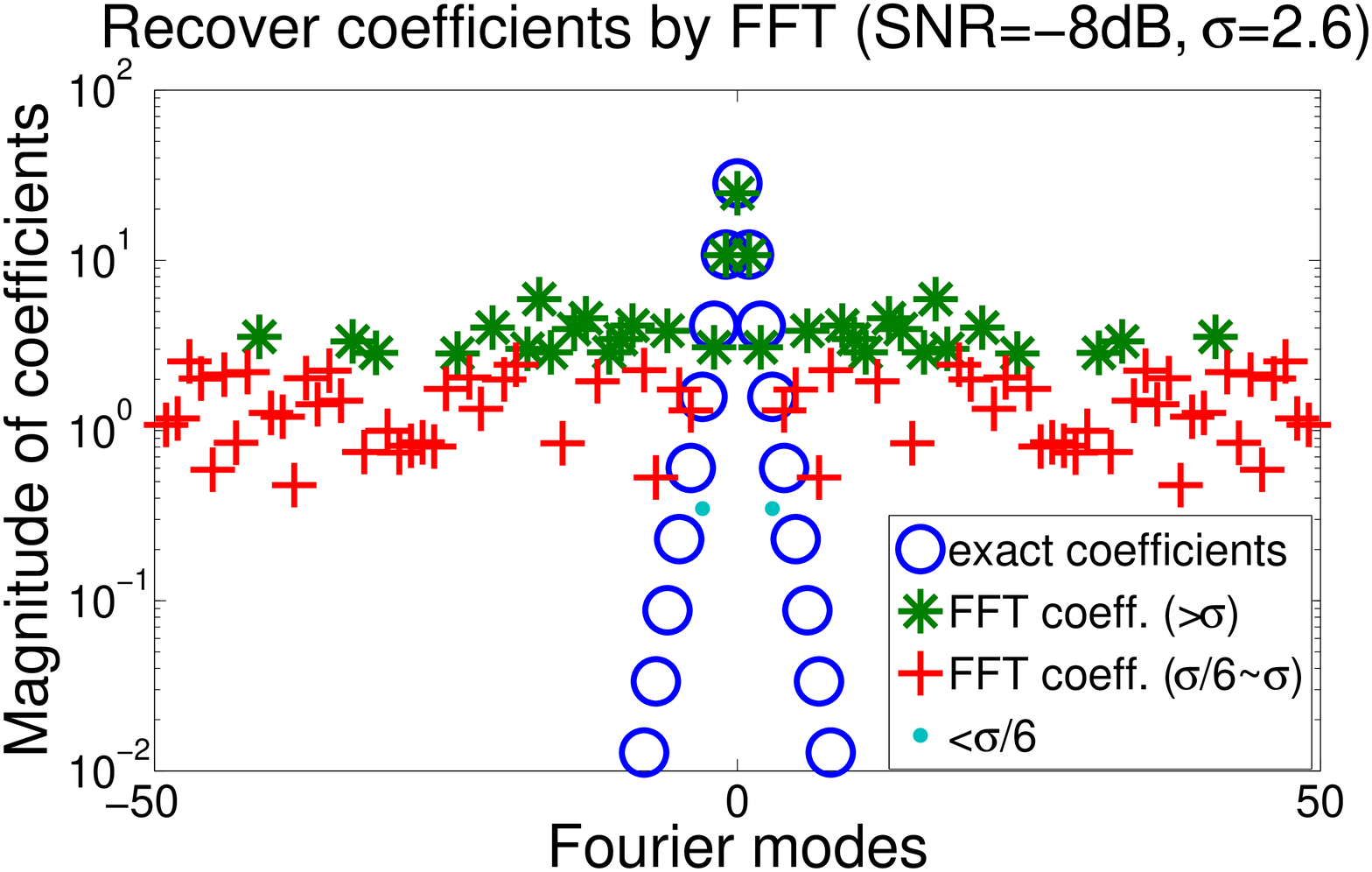}
\caption{For signal $S=1/(1.5+\cos 2\pi t)+noise$ with $SNR=-8dB$. Compare the approximation effect by RA$\ell$SFA and FFTW. Left: 
approximation of the significant coefficients 
by RA$\ell$SFA; the relative approximation error is $0.4\%$;right: approximation of the significant coefficients 
by FFTW. (this is for the one run illustrated. In other runs, it makes similarly one or two mistakes, not necessarily at the same modes.)} \label{fig:decay1}
\end{center}
\end{figure}    \addtocounter{table}{+1}

The results are shown in Figure \ref{fig:decay} ($SNR=-20dB$)and Figure \ref{fig:decay1} ($SNR=-8dB$), respectively.
For $SNR=-20dB$, the Fourier coefficients obtained by FFTW are all very close to the ``noise floor'', i.e., they
lie in a band of amplitude close to the value of $\sigma$. For $SNR=-8dB$, $\sigma$ is smaller ($\sigma=2.6$), and we find 
the ``noise floor'' in the FFTW computation at this lower level. The three largest modes of the signal have amplitudes significantly higher
than this $\sigma$, and FFTW finds them with reasonable accuracy. In contrast, RA$\ell$SFA (shown on the left
 in both figures; only 1 run is shown) hits all the coefficients exceeding $\sigma$ ``on the nose'', in both cases; it also finds all the central 15 modes exactly in the 
$SNR=-8dB$ case, even if they have values significantly smaller than $\sigma$. This experiment illustrates the great robustness
of RA$\ell$SFA  to noise and its ability to detect harmonic components with smaller energy than the white noise,
already seen in \ref{subset:noise}.

\subsection{Numerical Results in Two Dimensions}
\label{sect:ntwodim}
The number of grid points depends exponentially on the dimension. To achieve reasonable accuracy, a minimum $N$
is required in each dimension; however, when $d>1$, the FFTW has great difficulty in handling the corresponding $N^d$
points for even modest $N$. RA$\ell$SFA does not have this problem.

\subsubsection{Experiments for Eight-mode Signals in Two Dimensions}
We take the signal $S=\sum_{k=1}^{B}c_{k}\phi_{\omega_{x,k}}\phi_{\omega_{y,k}}$, where
$B=8,\epsilon=10^{-2}\|S\|,\delta=0.05$. The parameter $N$ is the number of grid points in
each dimension, random complex constants $c_k$ with real and imaginary parts in $[1, 10]$, and $\omega_{x,k}$ and $\omega_{y,k}$
are random integers from ${0, \ldots, N-1}$. As Table
\ref{tab:B2D2} shows, two dimensional RA$\ell$SFA surpasses two
dimensional FFTW when $N \geq 1500$. In particular, when $N=5000$ and the computation for samples is
 not included, the FFTW takes 21 seconds and RA$\ell$SFA only less than 5 second. When we include the
  sampling time, the crossover point becomes $N=900$.
The crossover point for $N$ is 70000 for $d=1$, and 900 for $d=2$; if we conjecture that the crossover $N$ for 2-mode in $d$ dimensions is given
by $c_2 n_2^{\frac{1}{d}}$, then this leads us to guess that the crossover $N$ for $d=3$ may be close to 210.

\begin{table}[htbp]
\begin{center}
\begin{tabular}{|c|c|c|c|c|c|c|}
\hline
Length& \multicolumn{2}{c|}{Time of}&  \multicolumn{1}{c|}{ Time of} & \multicolumn{2}{c|}{Time of RA$\ell$SFA}&\multicolumn{1}{c|}{Time of FFTW}\\
N &  \multicolumn{2}{c|}{RA$\ell$SFA}  &  \multicolumn{1}{c|}{FFTW} &\multicolumn{2}{c|}{(excluding sampling) }&
\multicolumn{1}{c|}{(excluding sampling)}\\ \cline{2-6}
    &clean& noisy& &clean& noisy& \\ \hline \hline
$100$&
 3.41& 3.64 & 0.05 &0.88 &1.05 &
0.04 \\ \hline $1000$&
 4.11 & 4.54 & 4.87  &1.04 &1.25 &
0.20 \\ \hline $2000$&
 4.76&
 4.91 & 20.86& 1.31 & 1.44&
2.12 \\ \hline $3000$&
 4.55 &
 5.37& 47.73& 1.33  & 1.70&
5.62\\ \hline $4000$&
 5.41 & 5.59&85.89 & 1.41&
 1.51 &
10.74\\ \hline $5000$&
 6.03 & 6.20& 138.27  & 1.56&
 1.66&
20.98\\ \hline
\end{tabular}
\end{center}
\caption{Time Comparison between RA$\ell$SFA and FFTW (B=8) based on 100 runs. ``Clean'' means that the test signal is pure. ``Noisy'' means the signal is contaminated with noise of $SNR=-4dB$. ``Excluding Sampling'' column lists the running time without including precomputation of sample values.} \label{tab:B2D2}
\end{table}   \addtocounter{figure}{+1}

%\begin{figure}[htbp]
%\begin{center}\includegraphics[%
%  width=6cm]{B82Dpureex.eps}
%\includegraphics[%
%  width=6cm]{quaB82Dpureex.eps} \\
%\end{center}
%\caption{Comparing the total running times of 8-mode 2-dimensional RA$\ell$SFA for 100 different runs of the randomized algorithm.
%Left: mean and variance as a function of $N$ (the number of points in each dimension); right: median, quantiles and total spread of the runs as a function
%of $N$, $B=8$, $d=2$} \label{fig:B2barD2}
%\end{figure}    \addtocounter{table}{+1}

\subsubsection{Experiments for Signals with Different Number of Modes $B$}
As in one dimension, the number of modes $B$ is the bottleneck for applying RA$\ell$SFA freely to signals that
are not so sparse. Suppose the signal is of the form $S(t)=\sum_{k=1}^{B}{c_{k}\phi_{\omega_{x,k}} \phi_{\omega_{y,k}}}$, with $N=3001$
for RA$\ell$SFA and $3000$ for FFTW. Table \ref{tab:diffBD2} illustrates the relationship between running time
and the number of modes $B$. Time increases depends polynomially on the number of terms $B$. When $N=3001$, the
crossover points for the FFTW to surpass RA$\ell$SFA are at $B=10$ and $B=17$ respectively, for including and
excluding sample computation cases. This implies the influence of $B$ on the execution time is far from
negligible.

\begin{table}[htbp]
\begin{center}
\begin{tabular}{|c|c|c|c|c|}
\hline
Number of modes & Time of & Time of & Time of RA$\ell$SFA &  Time of FFTW\\
B & RA$\ell$SFA & FFTW & (exclude sampling )& (exclude sampling)
\tabularnewline
\hline
$2$&
 0.15 &
 16.45 &
0.08& 5.64\tabularnewline $4$&
 0.52 &
 26.81 &
0.14& 5.64\tabularnewline $8$&
 4.55 &
 47.73 &
1.343& 5.64\tabularnewline $12$&
 19.37 &
 68.47&
8.82& 5.64\tabularnewline $16$&
 48.69 &
 89.13&
9.13 & 5.64\tabularnewline $20$&
 114.80 &
 109.88&
22.75& 5.64  \tabularnewline \hline
\end{tabular}
\end{center}
\caption{Time Comparison between RA$\ell$SFA and the FFTW for signals with different $B$ in 2 dimensions when
$N=3000$} \label{tab:diffBD2}
\end{table}   \addtocounter{figure}{+1}

\section{Theoretical Analysis and Techniques of RA$\ell$SFA}

We hope the numerical results have whetted the reader's appetite for
a more detailed explanation of the algorithm. Before explaining the structure of RA$\ell$SFA as implemented  by us, we review the basic idea of the algorithm.
Given a signal consisting of several frequency modes with different amplitudes, we could split it into several
pieces that have fewer modes. If one such piece had only a single mode, then it would be fairly easy to
identify this mode, and then to approximately find its amplitude. If the piece were not uni-modal, we could, by repeating
the splitting, eventually get uni-modular pieces. In order to compute the amplitudes, we need to ``estimate
coefficients.'' To verify the location of the modes in the frequency domain and concentrate on the most significant part of
the energy, we use ``group testing.'' An estimation that recurs over and over again
in this testing is the ``evaluation of norms.'' The first splitting of the signal is done in the ``isolation'' step.

The different steps are carried out on many different variants of the signals, each obtained by a random translation
in the frequency domain (corresponding to a modulation and the inverse dilation in the time
domain). Because the signal is sparse in the frequency domain, the different modes are highly likely to be well
separated after these random operations, facilitating isolation of individual modes.

The main skeleton of the algorithm was already given in \cite{GGIMS}; in our discussion here, we introduce new ideas and give the
corresponding theoretical analysis. We also explain how to set parameters that are either not mentioned or
loose in \cite{GGIMS}. In Section \ref{sect:num}, the total scheme of RA$\ell$SFA is given. In Section \ref{sect:coeff}, we show
the theoretical basis to choose parameters for estimating coefficients, and introduce some techniques to
speed up the algorithm. In Section \ref{sect:norm}, we set the parameters for norm estimation. Section \ref{sect:isolation} presents the
heuristic rules to pick the filter width for the isolation procedure. This is one of the key factors
determining the speed. A new filter is proposed for Group Testing in Section \ref{sect:group}, which works more efficiently.
Section \ref{sect:sample} discusses how to evaluate a random sample from a signal.
Finally, we discuss the extension to higher dimensions in Section \ref{sect:twodim}.

\subsection{Set-up of RA$\ell$SFA} \label{sect:num}

The following theorem is the main result of \cite{GGIMS}.

\begin{theorem}
Let an accuracy factor $\epsilon$, a failure probability $\delta$, and a sparsity target $B \in \mathbb{N}, B<N$
be given. Then for an arbitrary signal \(S\) of length $N$, RA$\ell$SFA will find, at a cost in time and space of
order \(poly(B,\log(N), \frac{1}{\epsilon}, \log(\frac{1}{\delta}))\) and with probability exceeding $1-\delta$, a
\(B-\)term approximation \(R\) to $S$, so that \(\|S-R\|^2 \leq (1+\epsilon)\|S-R_{opt}^B(S)\|^2_2\).
\end{theorem}

It is especially striking that the near-optimal representation $R$ can be built in sublinear time i.e.,
$poly(\log N)$ instead of the $O(N \log N)$ time requirement of the FFT. RA$\ell$SFA's speed will surpass the FFT as long as
the length of the signal is sufficiently large. In particular, if $S=R^B_{opt}(S)$ (that is, $\hat{S}(\omega)$
vanishes for all but $B$ values of $\omega$), then $\|S-R\|^2=0$, i.e., RA$\ell$SFA constructs $S$
without any error, at least in theory; in practice this means the error is limited by accuracy issues.

The main procedure is a Greedy Pursuit with the following steps:
\begin{algorithm} {\sc Total Scheme} \\
Input: signal $S$, the number of nonzero modes $B$ or its upper bound, accuracy factor $\epsilon$, success probability $1-\delta$, an upper bound of the signal energy $M$, the standard deviation of the white Gaussian noise $\sigma$, a ratio $\iota$ for relative precision.
\begin{enumerate}
    \item Initialize the representation signal $R$ to 0, set the maximum number of iterations $T=B\log(N)\log(\delta)/\epsilon^{2}$,
   \item Test whether $\Vert S-R \Vert^2 \leq \iota \|R\|^2$. If yes, return the representation signal $R$ and the whole algorithm ends; else go to step 3.
    \item Locate Fourier Modes $\omega$ for the signal $S-R$ by the isolation and group test procedures below.
    \item Estimate Fourier Coefficients at $\omega$: $\widehat{(S-R)}(\omega)$.
    \item If the total number of iterations is less than $T$, go to 2; else return the representation $R$.
  \end{enumerate}
\end{algorithm}

The test at stage 2, which is not in \cite{GGIMS}, can allow us to end early. The criterion $\|S-R\|^2 \leq \iota \|R\|^2$,
where $\iota$ is a small number chosen heuristically, is suitable when one expects that $S$ is sparse, up to a small energy contribution.
(Note that step 2 does not use the exact value of $\|S-R\|^2$, which is not known; we use
a procedure called norm estimation (see below) to give a rough estimate; this is good enough for the stop criterion. Other criteria could be substituted when appropriate.)

In practice, we would not know how many modes a signal has. In fact, the algorithm
doesn't really need to know $B$: it can just proceed until the residual energy
is estimated to be below threshold. (The value of $B$ is used only to set the
maximum number of iterations, and the width of a filter in the isolation
procedures below. For the maximum number $T$, a loose upper bound on $B$
suffices; the isolation filter width depends only very weakly on $B$.) If
either the residual energy or the threshold is large, the program would
continue. Note that for each iteration of the algorithm, we take new random
samples from the signal $S$.

\subsection{Estimate Individual Fourier Coefficients}
\label{sect:coeff}

The original RA$\ell$SFA only shows the validity of estimating coefficients without mentioning parameter
settings. Here we introduce a new technique to achieve better and faster estimation;
in the process, we give another proof
of Lemma 2 in \cite{GGIMS} that contains explicit parameter choices.

\begin{algorithm}\label{alg:coeff} {\sc Estimate Individual Fourier Coefficients} \\
Input: signal $S$, success probability $1-\delta$, and accuracy factor $\epsilon$.
\begin{enumerate}
    \item  Randomly sample from signal $S$ with indices
     $t_{i,j}$: $S(t_{i,j})$, $i=1,\ldots, 2\log(1/\delta)$, $j=0,\ldots,8/\epsilon^2$ .
    \item  Take the empirical mean of the $\left \langle S(t_{i,j}), \phi_{\omega}(t_j) \right \rangle$, $j=0,\ldots, 8/\epsilon^2$, store as $mean(i)$.
    \item  Take the median $y=median (mean(i))$, $i=1,\ldots, 2\log(1/\delta)$.
      \item Return $y$.
 \end{enumerate}
\end{algorithm} \label{alg:estcoef}

\begin{lemma} \label{lem:coef1}
Every application of Algorithm  {\em \ref{alg:coeff}} constructs a realization of a random variable $Z$, that
estimates the Fourier coefficient \(\hat{S}(\omega)\), good up to tolerance \(\epsilon^2 \|S\|^2\) with high
probability \(1-\delta\), i.e.,
\begin{equation}
Prob \left (|Z-\hat{S}(\omega)|^2 \geq \epsilon^2 \|S\|^2 \right ) \leq \delta.
\end{equation}
\end{lemma}
\begin{proof}
Define a random vector $V$ as follows:
\begin{equation}
V=\left (0,0,\cdots,NS(t),0,\cdots,0 \right )=N\delta_{t}S(t).
\end{equation}
where $t$ is chosen uniformly and randomly from $\{ l:\, l=1,\cdots,N\}$. Then the expectation of $V$ is
\begin{equation}
E(V)=\frac{1}{N}\sum_{t}{NS(t)\delta_{t}}.
\end{equation}
Let X be the random variable $X=\langle V,\phi_{\omega}\rangle$, where $\phi_{\omega}(t)=N^{-\frac{1}{2}}e^{-2\pi i\omega t/N}$. We have
\begin{equation}E[X]=\frac{1}{N}\sum_{t}{NS(t)\phi_{\omega}(t)}=\hat{S}(\omega),\end{equation} and
\begin{equation}E \left (|X-\hat{S}(\omega)|^{2}\right )\leq E(|X|^{2}) =\frac{1}{N}\sum_{t}{\left |\frac{N}{\sqrt{N}}S(t)e^{-2\pi i\omega t/N} \right |^{2}}=\Vert S\Vert_2^{2}.\end{equation}
Define another random vector $W$ as
the average of $L$ independent realization of $V$, with $L = 8 \epsilon^{-2}$. Let a random variable
\begin{equation}Y=\langle W,\phi_{\omega}\rangle.  \end{equation}
Then $E[Y]=\hat{S}(\omega)$ and $var[Y]=var[X]/L=\epsilon^{2}\Vert S\Vert^{2}/8$,
so that $Prob \left (|Y-\hat{S}(\omega)|^{2}\geq\epsilon^{2}\Vert S\Vert^{2} \right)\leq 1/8$. \\
Set $Z=median_K Y$, where $K=2\log(1/\delta)$. If $|Z-\hat{S}(\omega)|^2 \geq \epsilon^2 \|S\|^2$, then
for at least half of the $Y$s, we have
\begin{equation}
|Y-\hat{S}(\omega)|^2 \geq \epsilon^2 \|S\|^2.
\end{equation}
Therefore
\begin{align}
 P\left (|Z-\hat{S}(\omega)|^{2}\geq\epsilon^{2}\Vert S\Vert_2^{2} \right) &\leq \sum_{j=K/2} ^K {K \choose j} \left ( \frac{1}{8} \right )^j \nonumber \\
  &\leq 8^{-K/2} 2^K = 2^{-K/2} \leq\delta.
\end{align}
So with probability $1-\delta$, $Z$ is a good estimate of the Fourier Coefficient $\hat{S}(\omega)$, good up to
tolerance $\epsilon^{2}\Vert S\Vert^{2}.$  \qquad\end{proof}

Several observations and new techniques can speed up the coefficient estimation even further.

One observation is that fewer samples are already able to give an estimation with desirable accuracy and probability.
Our arguments indicate that $16 \epsilon^{-2} |\log (\delta)|$ samples per coefficient suffice to obtain good approximations of the
coefficients. The estimates used to obtain this bound are rather coarse, however. In a practical implementation, if a multi-step evaluation is used (see below),
it turns out that three steps, in which every step uses 10 samples per mean, and 5 means per median, for a total of 150 samples (per coefficient)
already determine the coefficient with accuracy $\epsilon = 10^{-4}$. The major factor in this drastic reduction (from $16\cdot 10^8|\log \delta|$ to 150)
is the much smaller number of means used; in practice, the dependence on $\epsilon$ grows much slower than $\epsilon^{-2}$ as $\epsilon \rightarrow 0$
If the signal is contaminated by noise or has more than one significant mode, we need more samples for a good estimation of the same accuracy.

An additional difference with the sampling described in \cite{GGIMS} is that one can replace individual random
samples by samples on short arithmetic progressions with random initial points. This technique became one of
several components in the RA$\ell$SFA version of \cite{GMS} that adapted the original algorithm in order to obtain
linearity in $B$. For a description of the arithmetic progression sampling, we refer to \cite{GMS}.
Surprisingly, this change not only improves the speed, but also gives a closer approximation than simply random
 sampling, using the same number of samples.

Another idea is a coarse-to-fine multi-step estimation of the coefficients. There are several reasons for not
estimating coefficients with high accuracy in only one step. One of them is that increasing the accuracy $\epsilon$
means a corresponding quadratic increase of the number of samples $O(|\log\delta| \epsilon^{-2})$.
A multi-step procedure, which produces only an approximate estimate of the coefficients in each step,
achieves better accuracy and speed. To explain how this works, we need the
following lemma.

\begin{lemma}
Given a signal $S$, let $\omega_1, \ldots, \omega_q$ be $q$ different frequencies, and define
$\beta \colon = \left [ \|S\|^2_2 - \sum_{i=1}^q |\hat{S}(\omega_i)|^2 \right ]/ \|S\|^2_2$.
Estimate the coefficients $\hat{S}(\omega_i)$ where
$i=1, \ldots, q$ by the following iterative algorithm: apply Algorithm {\em \ref{alg:coeff}} with precision $\hat{\epsilon}$ and probability of failure $\delta$;
keep the parameters fixed throughout the iterative procedure, and let $Z^n_i$, $i=1, \ldots, q$, be the estimate (at the $n$-th iteration)
of the $\omega_i$-th Fourier coefficient of $S-\sum_{k=1}^{n-1} \sum_{j=1}^q Z^k_j \phi_{\omega_j}$. The total estimate $R_n$ after the $n$-th
iteration is thus $R_n = \sum_{k=1}^n \sum_{j=1}^q Z^k_j \phi_{\omega_j}$. Then
\begin{equation}
\sum_{j=1}^q |\hat{S}(\omega_j) - \hat{R}_n(\omega_j)|^2 \leq \frac{q \hat{\epsilon}^2}{1-q \hat{\epsilon}^2} \beta \|S\|^2 + (q \hat{\epsilon}^2)^n \|S\|^2,
\end{equation}
with probability exceeding $(1-\delta)^{nq}$.
\end{lemma}
\begin{proof}
(This is essentially a simplified version of proof for Lemma 10 in \cite{GGIMS})\\
By Lemma \ref{lem:coef1},
\begin{equation}
|Z_i^n + \sum_{k=1}^{n-1} Z_i^k - \hat{S}(\omega_i)|^2 \leq \hat{\epsilon}^2 \|S-R_{n-1}\|^2,
\end{equation}
with probability exceeding $1-\delta$. It follows that
\begin{equation}
\sum_{i=1}^q |\hat{S}(\omega_i) - \sum_{k=1}^n Z^k_i|^2 \leq q \hat{\epsilon}^2 \|S-R_{n-1}\|^2,
\end{equation}
so that
\begin{align}
\|S-R_n\|^2 &\leq \sum_{\omega \notin \{ \omega_1, \ldots, \omega_q\}} |\hat{S}(\omega)|^2 + q \hat{\epsilon}^2 \|S-R_{n-1}\|^2 \nonumber \\
& = \|S\|^2 - \sum_{i=1}^q |\hat{S}(\omega_i)|^2 + q \hat{\epsilon}^2 \|S-R_{n-1}\|^2,  \\ \nonumber
& = \beta \|S\|^2 + q \hat{\epsilon}^2 \|S-R_{n-1}\|^2\label{eq:cof}
\end{align}
with probability exceeding $(1-\hat{\delta})^q$. \\
Consider now the sequence $(a_n)$, defined by $a_n = \beta \|S\|^2 + q \hat{\epsilon}^2 a_{n-1}$, where $a_0=\|S\|^2$. It is easy to see that
\begin{align}
a_n & = \beta \|S\|^2 \sum_{k=0}^{n-1} (q \hat{\epsilon}^2)^k + (q \hat{\epsilon}^2 )^n \|S\|^2 \\ \nonumber
& = \beta \|S\|^2 \frac{1-(q \hat{\epsilon}^2)^n}{1-(q \hat{\epsilon}^2)} + (q \hat{\epsilon}^2 )^n \|S\|^2.
\end{align}
It then follows by induction that $\|S-R_n\|^2 \leq a_n$, with probability exceeding $(1-\hat{\delta})^{nq}$, for all $n$; we have thus
\begin{align}
\|S-R_n\|^2 & \leq \beta \|S\|^2 \frac{1-(q \hat{\epsilon}^2)^n}{1-(q \hat{\epsilon}^2)} + (q \hat{\epsilon}^2 )^n \|S\|^2 \\ \nonumber
& \leq \beta \|S\|^2 \frac{1}{1-(q \hat{\epsilon}^2)} + (q \hat{\epsilon}^2 )^n \|S\|^2,
\end{align}
or equivalently,
\begin{equation}
\sum_{j=1}^q |\hat{S}(\omega_j) - \hat{R}_n(\omega_j)|^2 = \|S-R_n\|^2 - \beta \|S\|^2 \leq \beta \|S\|^2 \frac{q \hat{\epsilon}^2}{1-q \hat{\epsilon}^2}
+ (q \hat{\epsilon}^2 )^n \|S\|^2,
\end{equation}
with probability exceeding $(1-\delta)^{qn}$.  \qquad \end{proof}

The above lemma shows that repeated rough estimation can be more efficient than a single accurate estimation. To make this clear, if we set
\begin{equation}
q \epsilon^2 = \beta \frac{q \hat{\epsilon}^2 }{1-q \hat{\epsilon}^2 } + (q \hat{\epsilon}^2)^n, \, \, \, \,(1 - \delta)^q = (1-\hat{\delta})^{nq},
\label{eq:cofg}
\end{equation}
then a one-step procedure with parameters $\epsilon$, $\delta$ will achieve the same precision as an $n$-step iterative procedure with parameters
$\hat{\epsilon}$, $\hat{\delta}$. The one-step procedure will use $C q \epsilon^{-2} |\log(\delta)|$ sampling steps; the iterative procedure
will use $C nq \hat{\epsilon}^{-2} |\log(\hat{\delta})|$. It follows that the $n$-step iterative procedure will be more efficient, i.e., obtain the same
accuracy with the same probability while sampling {\it fewer} times, if
\begin{equation}
n \hat{\epsilon}^{-2} |\log (\hat{\delta})| \leq \epsilon^{-2} |\log (\delta)|,
\label{eq:cofgg}
\end{equation}
under the constraints (\ref{eq:cofg}). If $\beta=0$ (that is, if $S$ is a pure $q$-component signal), then this condition reduces 
(under the assumption that $\hat{\delta}$, $\delta$ and $\hat{\epsilon}$,
$\epsilon$ are small, so that $\frac{q \hat{\epsilon}^2 }{1-q \hat{\epsilon}^2 } \simeq q \hat{\epsilon}^2$, $(1-\hat{\delta})^n \simeq 1-n \hat{\delta}$) to
\begin{equation}
n \left ( |\log \delta| + n \right ) (q \hat{\epsilon}^2)^{n-1} \leq |\log \delta|,
\label{eq:cofg1}
\end{equation}
which is certainly satisfied if $\hat{\epsilon}$ is sufficiently small and $n$ sufficiently large. If $\beta \neq 0$, matters are more complicated,
but by a simple continuity argument we expect the condition still to be satisfied if $\beta$ is sufficiently small. If $\beta$
is too large, (e.g. if $\beta > n_0^{-1}$, where $n_0$ is the minimum value of $n$ for which (\ref{eq:cofg1}) holds), then there are no choices of $n$,
$\hat{\epsilon}$, $\hat{\delta}$ that will satisfy (\ref{eq:cofg}) and (\ref{eq:cofgg}). On the other hand, $\beta$ can be large only if $S$ has important modes not included in ${\omega_1, \ldots, \omega_q}$. In practice, we use the multi-step procedure after the most important modes have been identified so that $\beta$ is small. For sufficiently small $\beta$, we do gain by taking the iterative
procedure. For example, assume that $\beta=10^{-2}$, for a signal of type $S = \phi_1+\phi_2$ with $N=1000$, $q=B=2$, $\delta=2^{-7}$,
$\epsilon=4 \cdot 10^{-4}$,
and with $n=3$, theoretically we would then use 450,000 samplings for the one-step
procedure, versus 150 samples for the iterative procedure. Note that we introduced the parameter $\beta$ only for expository
purposes. In practice, we simply continue with the process of identifying modes and roughly estimating their coefficients until our estimate of the residual signal is small; at that point, we switch to the above multi-step estimation  procedure.

\subsection{Estimate Norms}
\label{sect:norm}

The basic principle to locate the label of the significant frequency is to estimate the energy of the new
signals obtained from isolation and group testing steps. The new signals are supported on only a small number of taps in
the time domain and have 98$\%$ of their energies concentrated on one mode. The original analysis in
\cite{GGIMS} only gave its loose theoretical bound. Here we find the empirical parameters, i.e., the number of
samples for norm estimation.

Here is a new scheme for estimating norms, which uses much fewer samples than the original one and still achieves good estimation.
It can ultimately be used to find the significant mode in conjunction with Group Testing and MSB, below.

\begin{algorithm} \label{alg:norm} {\sc Estimate Norms}
Input: signal $S$, failure probability $\delta$.
\begin{enumerate}
    \item  Initialize: the number of samples: $r=\lfloor 12.5 \ln (1/\delta) \rfloor $.
    \item  Take $r$ independent random samples from the signal $S$: $S(i_1), \ldots, S(i_r)$, where $r$ is a multiple of 5.
    \item  Return $ N \times$ ``60-th percentile of'' ${|S(i_1)|^2, \ldots, |S(i_r)|^2}$.
 \end{enumerate}
\end{algorithm}

The following lemma presents the theoretical analysis of this algorithm.

\begin{lemma}\label{lem:norm}
If a signal $S$ is $93\%$ pure, the number of samples $r>12.5 \ln (1/\delta)$, the output of Algorithm {\em \ref{alg:norm}} gives an estimation $X$ of its energy
which exceeds $0.3 \|S\|^2$ with probability exceeding $1-\delta$.
\end{lemma}
\begin{proof}
Without loss of generality, suppose that $\|S\|=1$. Suppose the signal $S=a\phi_{\omega}+e$, where $|a|^2>0.93 \|S\|^2$, and $\phi_{omega}$ and $e$ are orthogonal. 
We shall sample the signal $S$ independently for $r$ times, as stated in Algorithm \ref{alg:norm}. Note that we do not impose that
samples be taken at different time positions; with very small probability, the samples could coincide. Let $T=\{t:N|S(t)|^2< 0.3 \|S\|^2 \}$.
Hence, for any $t\in T$, we have $\sqrt{N} |S(t)|<\sqrt{0.3}=0.5477$. Also by the purity of $S$, we have $\|e\|^2 \leq 0.07$. Since $|S(t)|\geq |a \phi_{\omega}(t)|-|e(t)|$, we obtain
\begin{equation}
\sqrt{N}|e(t)|>|a|-\sqrt{N}|S(t)|.
\end{equation}
then for any $t\in T$,
\begin{equation}
\sqrt{N}|e(t)|>\sqrt{0.93}-\sqrt{0.3}.
\end{equation}
Therefore,
\begin{equation}
0.07N \geq N \|e\|^2 \geq N \sum_{t\in T}|e(t)|^2 \geq (\sqrt{0.93}-\sqrt{0.3})^2 |T|.
\end{equation}
It follows that
\begin{equation}
|T| \leq 0.403 N
\end{equation}
Let $\alpha = \frac{|T|}{N}$; the above inequality becomes
$0 \leq \alpha \leq 0.403$.  \\
Consider now the characteristic function $\chi_T$ of the set $T$, 
\begin{equation}
\chi_T(t) = \begin{cases}
1 & \text{if $t \in T$} \\
0 & \text{otherwise},
\end{cases}
\end{equation}
and define the random variable $X_T$ as $\chi_T(i)$, where $i$ is picked randomly. Then we have
\begin{equation}
E(X_T)=\frac{|T|}{N} \leq 0.403,
\end{equation}
and
\begin{equation}
E(e^{X_T z}) = e^0 Prob(\chi_T(i)=0) + e^z Prob(\chi_T(i)=1) = 1-\alpha + \alpha e^z.
\end{equation}
Suppose now we sample the signal $S$ $r$ times independently, and obtain $S(t_1), \ldots, S(t_r)$, where $t_1, \ldots, t_r \in [0,N]$. Take the \mbox{\it {60-th}} percentile of the numbers $N|S(t_1)|^2, \ldots, N|S(t_r)|^2$.
By Chernoff's standard argument, we have for $z>0$
\begin{align}
Prob \left ( \text{60-th percentile} < 0.3 \|S\|^2 \right ) & = Prob \left (0.6r\,\, \text{of the samples'} \,\,
t \,\, \text{belong to T}  \right ) \nonumber \\
& = Prob(\chi_T(t_1)+ \ldots + \chi_T(t_r) > 0.6r) \nonumber \\
& \leq e^{-0.6rz} E(e^{z \sum_{j=1}^r \chi_T(t_j)} ) \nonumber \\
&= \left [ (1-\alpha) e^{-0.6z} + \alpha e^{0.4 z} \right ]^r.
\end{align}
Take $z=\ln (1.5(1-\alpha)/ \alpha)$, then
\begin{equation} \label{eq:norm}
(1-\alpha)e^{-0.6z} + \alpha e^{0.4z} = 1.96 \alpha^{0.6} (1-\alpha)^{0.4}.
\end{equation}
The right hand side of (\ref{eq:norm}) is increasing in $\alpha$ on the interval $[0, 0.403]$; since $\alpha \leq 0.403$, 
we obtain an upper bound by substituting 0.403 for $\alpha$:\begin{eqnarray}
\left [ (1-\alpha) e^{-0.6z} + \alpha e^{0.4 z} \right ]^r = \left [ 1.96 \alpha^{0.6} (1-\alpha)^{0.4} \right ]^r \leq e^{-0.08r}.
\end{eqnarray}
So for $r \geq 12.5 \ln (1/\delta)$, we have
\begin{align}
Prob(\text{Output of Algorithm}\, 3.6 \geq 0.3 \|S\|^2)& =Prob( \text{60-th percentile of}\, N|S(t)|^2 \geq 0.3 \|S\|^2) \\ \nonumber
& \geq 1-\delta.
\end{align} \qquad\end{proof}

In practice, we often generate signals that are not so pure and thus need more samples for norm estimation.
Although the estimation is sometimes pretty far away from the true value, it gives a rough idea of where the significant mode might be. When we desire more accuracy,
a smaller constant $C$ in the number of samples $C\log(1/\delta)$ is chosen. In the statement of the algorithm, we choose $r$ to be a multiple of 5, so that
the \mbox{\it {60-th}} percentile would be well-defined. In practice, it works equally well to take $r$ that are not multiples
of 5 and to round down, taking the $\lfloor 3r/5 \rfloor$-th sample in an increasingly ordered set of samples.

We shall also need an upper bound on the outcome of Algorithm \ref{alg:norm}, which should hold regardless of whether
the signal $S$ is highly pure or not. This is provided by the next lemma, which proves that for general signals, Algorithm \ref{alg:norm} produces an estimation of the energy, that is less than $2\|S\|^2$ with high probability.

\begin{lemma} \label{lem:norm2}
Suppose Algorithm \ref{alg:norm} generates an estimation $X$ for $\|S\|^2$, then 
\begin{equation}
Prob(X \geq 2\|S\|^2) \leq  \left ( \frac{1 }{2} \right )^{0.144 \ln(1/\delta)} = \delta^{0.1}.
\end{equation}
\end{lemma}
\begin{proof}
Suppose $r$ independent random samples are $S(t_1), S(t_2), \ldots, S(t_r)$, then
\begin{equation}
Prob( N|S(t_i)|^2 \geq 2 \|S\|^2) \leq \frac{ N E(|S(t_i)|^2) }{ 2 \|S\|^2} = 1/2.
\end{equation}
Since $X$ is the 60-th percentile of the sequence $NS(t_1), \ldots, NS(t_r)$, with $r=0.36 \ln (1/\delta)$, 
 \begin{equation}
Prob( X \geq 2 \|S\|^2) \leq \left( Prob( N|S(t_i)|^2 \geq 2 \|S\|^2) \right)^{0.144  \ln(1/\delta)} \leq \left ( \frac{1 }{2} \right )^{0.144 \ln(1/\delta)} = \delta^{0.1}.
\end{equation}
 \qquad\end{proof}

\subsection{Isolation}
\label{sect:isolation}

Isolation processes a signal $S$ and returns a new signal with significant frequency $\omega$, with 98$\%$ of
the energy concentrated on this mode. A frequency $\omega$ is called ``significant'' for $S$ , if
$|\hat{S}(\omega)|>\eta\|S\|^2$, where $\eta$ is a threshold, fixed by the implementation, which may be fairly
small. More precisely, the isolation step returns a series of signals $F_0, F_1, \ldots, F_{r}$, such that, with high probability,
$|\hat{F}_j(\omega)|^2 \geq 0.98 \|F_j\|^2$ for some $j$, that is, at least one of the $F_0, F_1, \ldots, F_{r}$ is $98\%$ pure.

Typically, not all of the $F_i$s are pure. We shall nevertheless apply the further steps of the algorithm to each of the
$F_i$s, since we don't know which one is pure. An impure $F_i$ may lead to a meaningless value for the putative mode
$\tilde{\omega}_i$ located in $F_i$. This is detected by the computation of the corresponding coefficients: only when
the coefficient corresponding to a mode is significant do we output the mode and its coefficient. Some impure signals might output an insignificant mode.
Hence, we estimate and compare their coefficients to check the significance of the modes.
Finally, we only output the modes with significant coefficients.

The discussion in \cite{GGIMS} proposes a B-tap box-car filter in the time domain, which corresponds to a
Dirichlet filter with width $\frac{N}{B}$ in the frequency domain. The whole frequency region would be covered
by random dilation and translations of this filter.

Notation: as in \cite{weaver}, we define a box-car filter $H_k$ as
$H_{k}(t)=\frac{\sqrt{N}}{2k+1}\chi_{[-k,k]}$, where $k \in \mathbb{N}$.

\begin{lemma}
\begin{enumerate}
    \item For all \(k\),
                 \begin{equation}\label{ji}
                  \hat{H}_{k}(\omega) = \frac{1}{2k+1}\sum_{t=-k}^{k}e^{\frac{-2 \pi i \omega t}{N}}
                   = \frac{\sin(\pi(2k+1)\omega/N)}{(2k+1)\sin(\pi \omega/N)}.
               \end{equation}
  \item Notation: $H_{k,j}(t) = e^{2 \pi i j t /(2k+1)}H_k(t)$ in the time domain, which is equivalent to a shift of $\hat{H}_k(\omega)$ by $jN/(2k+1)$ in the
frequency domain.
\item Notation: Define $R_{\theta, \sigma}S(t)$ by $R_{\theta, \sigma}S(t)= e^{-2 \pi i \theta t/\sigma N} F(t/\sigma)$, so that 
$\widehat{R_{\theta, \sigma} S}=\hat{S}(\sigma \omega + \theta)$., where $\widehat{R_{\theta,
\sigma}}$ is a dilation and shift operator in the frequency domain.
\end{enumerate}
\end{lemma}

More detailed description of the Box-car filter can be found in \cite{GGIMS}.

The isolation procedure in \cite{GGIMS} randomly permutes the signal $S$ and then convolves it with a shifted
version of $H_{k,j}$ to get a series of new signals $F_j = H_{k,j}* R_{\theta, \sigma} S$, where $j=0, \ldots, 2k$.
This scheme does not work well in practice. In the new version of the isolation steps, each $F_j = H_k
* R_{\theta_j, \sigma_j} S$ corresponds to different randomly generated dilation and modulation factors,
with $j=0, \ldots, \log(1/\delta)$, the parameters $\sigma_j$ and $N$ are relatively prime. These factors are taken at random between 0 and $N-1$. The following lemma is similar to Lemma 8 in \cite{GGIMS} for
the new isolation step, with more explicit values of the parameters.

\begin{lemma} {\em \cite{GGIMS}}
Let a signal $S$ and a number $\eta$ be given, and create $\log(1/\delta)$ new signals: $F_0, \ldots,
F_{\log(1/\delta)}$ with $F_j = H_k*R_{\theta_j, \sigma_j}S$, where $j=0, \ldots, \log (1/\delta)$. If $k \geq 12.25 (1-\eta)\pi^2/\eta$
, then for each $\omega$ such that $|\hat{S}(\omega)|^2 \geq \eta \|S\|^2$, there
exists some $j$ such that  with high probability $1-\delta$, the new signal $F_j$ is $98\%$ pure.
\end{lemma}
\begin{proof}
Suppose $\sigma_j^{-1} (\omega - \theta_j)$ falls into the pass region of the $H_k$ filter, i.e.,  that $\left |\sigma_j^{-1} (\omega -
\theta_j \right | \leq \frac{N}{2(2k+1)}$.
We know that
\begin{equation}
\left | \hat{H}_k \left (\sigma_j^{-1} (\omega -
\theta_j) \right ) \right | \geq 2/\pi,
\end{equation}
so that
\begin{equation}
\left |\hat{F}_j \left (\sigma_j ^{-1} (\omega-\theta_j) \right ) \right |^2 \geq (2/\pi)^2 \left
|\hat{S}(\omega) \right |^2 \geq (2/\pi)^2 \eta \|S\|^2.
\end{equation}
greater than the average value, $1/(2k+1)$, of $|H^k|^2$. 
Since $|\hat{H}_k(\sigma_j^{-1}(\omega-\theta_j))|^2$ is greater than the average value of $\hat{H}_k$, we have
\begin{equation}
\frac{\sum_{\omega' \neq \sigma_j^{-1}(\omega-\theta_j)} |\hat{H}_k(\omega')|^2}{N-1} \leq \frac{\|H_k\|^2}{N}=
\frac{1}{2k+1}.
\end{equation}
Moreover, $\sum_{\omega'' \neq \omega} |\hat{S}(\omega')|^2 \leq (1-\eta)\|S\|^2$. In
particular, $|\hat{S}(\omega')|^2 \leq (1-\eta)\|S\|^2$ if $\omega' \neq \omega$. We then have
\begin{equation}
E\left[\sum_{\omega'\neq\sigma_j^{-1}(\omega-\theta_j)}  | \hat{F}_j(\omega')|^{2}  \bigg |
-\frac{1}{2}N/(2k+1)\leq\sigma_j^{-1}(\omega-\theta_j)\leq \frac{1}{2}N/(2k+1)\right]\leq\frac{(1-\eta)\|
S\|^{2}}{2k+1}.
\end{equation}
Define $X$ to be the random variable
\begin{equation}
X= \left  \{ \sum_{\omega' \neq \sigma_j^{-1}(\omega-\theta_j)}  | \hat{F}_j(\omega')|^{2}  \bigg |
-\frac{1}{2}N/(2k+1)\leq\sigma_j^{-1}(\omega-\theta_j)\leq \frac{1}{2}N/(2k+1) \right \}.
\end{equation}
 For this random variable, we have
\begin{eqnarray}
Prob \left ( \frac{X}{|\hat{F}_{j}(\sigma_j^{-1}(\omega-\theta_j))|^{2}}\geq 1/49 \right ) &=Prob \left (X\geq|\hat{F}_{j}(\sigma_j^{-1}(\omega-\theta_j))|^{2}/49\right ) \nonumber \\
 & \leq\frac{E(X)}{|\hat{F}_{j}(\sigma_j^{-1}(\omega-\theta_j))|^{2}/49}\leq\frac{49(1-\eta)\pi^{2}}{4 \eta (2k+1)}.
\end{eqnarray}
Since $k \geq 12.25 (1-\eta)\pi^2/\eta$, the right hand side of (4.37) is $\leq 1/2$, meaning that the signal $F_j$  is 98$\%$ pure with
probability $\geq 1/2$.
 The success probability, i.e.,  the probability of obtaining at least one $F_j$ that is 98$\%$ pure, can be boosted from $\frac{1}{2}$
to probability $1-\delta$ by repeating $O(\log(1/\delta))$ times, i.e.,  generating $O(\log(1/\delta))$ signals.
 \qquad\end{proof}

The above lemma gives a lower bound for the filter width. Obviously, the larger the width in the time domain,
the higher the probability that the frequency will be successfully isolated. However, a larger width leads to
more evaluations of the function and therefore more time for each isolation step. One needs to balance carefully
between the computational time for each iteration step and the total number of iterations.

Based on several numerical experiments, we found that a very narrow filter is preferable and gives good
performance; for instance, the filter with three-tap width, i.e.,  $k=1$
works best for a signal with 2 modes. For the choice $k=4$, the algorithm ends after fewer iterations; however,
each iteration takes much more time. The choice of a 9-tap width filter makes the code four times
slower in total.

The filter width is weakly determined by the number of modes in the signal, not by the length of the signal.
Through experimentation, we found that when the number of modes is less than 8, the 3-tap width filter works
very well; as the number of modes increases, larger width filters are better. Numerical experiments suggests
a sublinear relationship between the width of the filter and the number of modes; in our experiments a 5-tap filter still
sufficed for $B=64$.

\subsection{Group Testing}
\label{sect:group}
After the isolation returns several signals, at least one of which is 98$\%$ pure with high probability, group
testing aims at finding the most significant mode for each. We use a procedure called Most
Significant Bit (MSB) to approach the mode recursively.

In each MSB step, we use a Box-car filter $H_k$ to subdivide the whole region into $2k+1$ subregions. By
estimating the energies and comparing the estimates for all these new signals, we find the one with maximum
energy, and we exclude those that have estimated energies much smaller than this maximum energy. We then repeat
on the remaining region, a more precisely on the region obtained by removing the largest chain of excluded
intervals; we dilate so that this new region fills the whole original interval, and split again. The successive
outputs of the retained region gives an increasingly good approximation to the dominant frequencies. The following are the Group testing procedures:

\begin{algorithm} \label{alg:group}{\sc Group Testing} \\
Input: signal $F$, the length $N$ of the signal $F$. \\
Initialize: set the signal $F$ to $F_0$, iterative step $i=0$, the length $N$ of the signal, the accumulation factor $q=1$. \\
In the $i$th iteration,
 \begin{enumerate}
    \item  If $q \geq N $, then return 0.
    \item  Find the most significant bit $v$ and the number of significant intervals $c$ by the procedure MSB.
    \item  Update $i=i+1$, modulate the signal $F_i$ by $\frac{(v+0.5)N}{4(2k+1)}$ and dilate it by a factor of $4(2k+1)/c$. Store it in $F_{i+1}$.
    \item  Call the Group testing again with the new signal $F_{i}$, store its result in $g$.
    \item  Update the accumulation factor $q = q * 4(2k+1)/c$.
      \item If $g> N/2$, then $g = g -N$.
      \item return $mod( \lfloor \frac{cg}{4(2k+1)}+ \frac{(v+1/2)N}{4(2k+1)}+0.5 \rfloor, N)$;
     \end{enumerate}
\end{algorithm}

The MSB procedure is as follows.
\begin{algorithm} \label{alg:msb}{\sc MSB (Most Significant Bit)} \\
\text{}\hspace{10mm} Input: signal $F$ with length $N$, a threshold $0<\eta<1$.
\begin{enumerate}
    \item Get a series of new signals $G_j(t) =F(t) \star (e^{2 \pi i j t/4(2k+1)} H_k )$, $j=0, \ldots, 8k+4$. That is, each signal $G_j$ concentrates on the
pass region $[ \frac{(j-1/2)N}{4(2k+1)}, \frac{(j+1/2)N}{4(2k+1)}]:=pass_j$.
    \item Estimate the energies $e_j$ of $G_j$, $j=0, \ldots, 8k+4$.
    \item Let $l$ be the index for the signal with the maximum energy.
      \item Compare the energies of all other signals with the $l$th signal. If $e_i < \eta e_l$
, label it as an interval with small energy.
    \item Take the center $v_s$ of the longest chain of  consecutive small energy intervals, suppose there are $c_s$ intervals altogether in this chain.
      \item The center of the large energy intervals is $v = 4(2k+1)-v_s$, the number of intervals with large energy is $c = 4(2k+1)-c_s$.
    \item If $c>4(2k+1)/2$, then do the original MSB {\em \cite{GGIMS}} to get $v$ and set $c=2$, and $v=
    center\, of\, the\, interval\, with\, maximal\, energy$.
    \item Output the dilation factor $c$ and the most significant bit $v$.
\end{enumerate}
\end{algorithm}

\begin{lemma}
Given a signal \(F\) with \(98\)$\%$ purity, suppose \(G_j(t) =F * e^{2 \pi i j t /4(2k+1)} H_k(t)\). If \(k \geq 2
\), then Algorithm {\em \ref{alg:group}} can find the significant frequency $\omega$ of the signal $F$ with high
probability.
\end{lemma}
\begin{proof}
Suppose the filter width of $H_k$ is $2k+1$. Observe that, for some $j$, $0\leq j \leq 4(2k+1)$, $\omega \in pass_j$.
Without loss of generality, assume $j=0$. Now consider the signal $G_0$. Since $\omega \in pass_0$, the Fourier coefficient 
$\hat{G}_0(\omega)$ satisfies
\begin{align}
|\hat{G_0}(\omega)|^2 & \geq \left( \frac{\sin(\pi/8)}{(2k+1)sin(\pi/8(2k+1))} \right ) ^2  |\hat{F}(\omega)|^2  \\ \nonumber
 & \geq \left( \frac{\sin(\pi/8)}{(2k+1)\sin(\pi/8(2k+1))}\right)^2 (0.98) \|F\|^2 \\ \nonumber
& \geq 0.9744^2 \cdot 0.98 \|F\|^2 \approx 0.93 \|F\|^2.
\end{align}
for all $k>0$. It follows from Lemma \ref{lem:norm}, that the output of Algorithm \ref{alg:norm}, applies to $G_0$, estimate that is at least
\begin{equation}
0.3\|G_0\|^2 \geq 0.3 |\hat{G_0}(\omega)|^2 \geq  0.3 \cdot 0.98 \left( \frac{\sin(\pi/8)}{(2k+1)\sin(\pi/8(2k+1))}\right)^2 \|F\|^2 .
\end{equation}
On the other hand, now consider $G_5$. Note that
\begin{align}
|\hat{G_5}(\omega)| = |\hat{F}(\omega)||\widehat{H_{k}}(\omega)| & \leq  \frac{1}{(2k+1)\sin(9 \pi/8(2k+1))} |\hat{F}(\omega)| \\ \nonumber
  &\leq  \frac{1}{(2k+1)\sin(9\pi/8(2k+1))} \|F\|.
\end{align}
Also,
$\| G_{5}\|^{2}-|\hat{G}_{5}(\omega)|^{2} \leq 0.02\| F\|^{2}$, because $F$ is $98\%$ pure. Thus 
\begin{equation}
\| G_{5}\|^{2}\leq|\hat{G}_{5}(\omega)|^{2}+0.02\| F\|^{2}=\left( \left( \frac{1}{(2k+1)\sin(9 \pi/8(2k+1))}\right)^2 +0.02  \right) \| F\|^{2}.
\end{equation}
By Lemma \ref{lem:norm2}, if we use Algorithm \ref{alg:norm}, the estimation result for $G_5$ will be at most $2 \| G_{5}\|^{2}$ with high probability.
It is easy to show that the inequality
\begin{equation} 
0.294   \left( \frac{\sin(\pi/8)}{(2k+1)\sin(\pi/8(2k+1))}\right)^2 \geq 2 \left( \frac{1}{(2k+1)\sin(9 \pi/8(2k+1))}\right)^2 +0.04 
\end{equation}
holds for all $k>0$. The same argument applies to $G_j$ with $5\leq j\leq 4(2k+1)-5$. It follows that, with high probability, the result of applying Algorithm 
\ref{alg:norm} to $G_0$ will give a result that exceeds the result obtained by applying Algorithm \ref{alg:norm} to $G_j$ with $5 \leq j \leq 4(2k+1)-5$.

In general, if the pass region is at some $j_{0}$, we can compare $\| G_{j_{0}}\|^{2}$with
$\| G_{j}\|^{2}$for all $|j-j_{0}|\geq5$. If there is some $j_{0}$ for which
the estimation of $\| G_{j_{0}}\|^{2}$ is apparently larger than
$\| G_{j}\|^{2}$, then we conclude $\omega\notin pass_{j}$; otherwise, possibly
$\omega\in pass_{j}$. By the above argument, we can eliminate
$4(2k+1)-9$ consecutive pass regions out of the $4(2k+1)$, leaving a cyclic interval
of length at most $\frac{9N}{4(2k+1)}$. In order for the residual region to be smaller or equal 
to half of the whole region, we need $4(2k+1) \geq 18$, which is equivalent to the condition $k\geq 2$.

In the recursive steps, let $P$ denote a cyclic interval with size
at most $\frac{9N}{4(2k+1)}$ that includes all the possibilities for $\omega$.
Let $v$ denote its center. Then generate a new signal $F_{1}(t)=e^{-2\pi ivt/N}F(t)$; this
is a shift of the spectrum of $F$ by $-v$. Thus the frequency $\omega-v$
is the biggest frequency of $F_{1}(t)$., which is in the range of
$-\frac{4.5N}{4(2k+1))}$ to $+\frac{4.5N}{4(2k+1))}$. We will now seek $\omega-v$.

Since we rule out a fraction of $\frac{(8k-5)N}{4(2k+1)}$ length of the whole region, we may
dilate the remainder by $\lfloor 4(2k+1)/9 \rfloor$, which can be accomplished in
the time domain by dilating $F_{1}$ by $\frac{9}{4(2k+1)}$. Thus the interval of
length just less than $\frac{9N}{4(2k+1)}$ known to contain $\omega-v$ is dilated
to the alternate positions in an interval of length just less than
$N$. We then rule out again $\frac{8k-5}{4(2k+1)}$ of this dilated frequency
domain, leaving a remainder of length at most $\frac{9}{4(2k+1)}$ length. Then we undo the dilation,
getting an interval of length just less than $\frac{9N}{(4(2k+1))}$, centered
at some $v_{2}$, which is the second most significant bit of $\omega$ in
a number base $\lfloor \frac{4(2k+1)}{9} \rfloor$. We would repeat this process to get the other bits
of $\omega$. By getting a series of $v_{1},\ldots,v_{\lfloor\log_{4(2k+1)/9}N\rfloor+1}$,
we can recover the $\omega$. \qquad\end{proof}

In fact, a narrower filter with a larger shift width than $\frac{N}{4(2k+1)}$ works fine and makes the algorithm faster in practice.
 Heuristically, we find that the optimal number of
taps for small $B$ cases is 3. Suppose the MSB filter width is 3 and each MSB rules out 2 intervals out of 3,
then the total number of recursive group test is $\log_3 N$. Then the computational cost is $3 \log_3 N$ norm estimations and $2 \log_3 N$ comparisons. Numerical experiments suggests that $k$ is probably linear in $\log B$. The shift width we use in practice is $\frac{N}{2k+1}$. 

We find that the output of group testing in both the original and the present version of RA$\ell$SFA might differ from the true mode by one place.
We suspect that the reason is that all the float operations and the conversion to integers introduce and accumulate some error into the final
frequency. As a solution, the coefficients of nearby neighbors are also estimated roughly to determine the true significant modes.

\subsection{Sample from a transformed signal}
\label{sect:sample}

A key issue in the implementation consists of
obtaining information (by sampling) from a signal after it has been dilated,
modulated, or even convolved. We briefly discuss here how to carry out this
sampling in discrete signals.

First, we consider a dilated and modulated signal,
for example, in the isolation procedure which uses $R_{\theta,
\sigma}S(t)=e^{-2\pi i \theta t/\sigma N}S(t/\sigma)$, which is equivalent to
$\widehat{(R_{\theta, \sigma}S)}(\omega) = \hat{S}(\sigma \omega + \theta)$ in
the frequency domain. Here $\sigma$ and $\theta$ are chosen uniformly and
randomly, from $0$ to $N-1$ for $\theta$, and from $1$ to $N-1$ for $\sigma$.
The sample $R_{\theta, \sigma}F(t)$, where $t \in \lbrace 0, 1, \ldots,
N-1\rbrace$,is then $e^{-2 \pi i \theta t /\sigma N}(R_{\theta, \sigma})F(t)=e^{-2 \pi i
\theta t /\sigma N}F(\sigma^{*}t)$, where $\sigma^*$ is chosen so that $\sigma^{*}
\sigma =1(mod\, N)$. If $N$ is prime, then we can always find (a unique value
for) $\sigma^{*}$ for arbitrary $\sigma$; if $N$ is not prime, $\sigma^{*}$ may
fail to exist for some choices of $\sigma$. Our program uses the Euclidean
algorithm to determine $\sigma^{*}$; when $N$ is not prime and $\sigma$ and $N$
are not co-prime, the resulting candidates for $\sigma^{*}$ are not correct and
may lead to estimates for the modes that are incorrect; these mistakes are
detected automatically by the algorithm when it estimates the corresponding
coefficient and finds it to be below threshold.

We also need to sample from convolved signals, e.g.  $S*H_k(t)$. Because $H_k$ has only $2k+1$ taps, only
$2k+1$ points contribute to the calculation of the convolution. Since
$S*H_k(t)=\sum_{i=-k}^k H_k(i) S(t-i)$, we need only the values $S(t-i)$,
$i=-k,\ldots,k$, all of which we sample.

\subsection{Extension to a Higher Dimensional Signal}
\label{sect:twodim}

The original RA$\ell$SFA discusses only the one dimensional case. As explained earlier, it is of particular interest
to extend RA$\ell$SFA to higher dimensional cases because there its advantage over the FFT is more pronounced.

In $d$ dimensions, the Fourier basis function is
\begin{equation}
\phi_{\vec \omega (\vec x) } = \phi_{\omega_1, \ldots, \omega_d}(x_1, \ldots, x_d) = N^{-\frac{d}{2}} e^{i 2 \pi
\omega_1 x_1 / N +\ldots +i 2 \pi \omega_d x_d / N } = N^{-\frac{d}{2}} e^{i 2 \pi \vec \omega_i \vec x_i /N};
\end{equation}
 the representation of a signal is
 \begin{equation}
 S(x_1, \ldots, x_d) = \sum_{i=1} ^N {c_i \phi_{\omega_{i,1},\ldots,\omega_{i,d}}}.
 \end{equation}
  Suppose the dimension of the signal is $d$, denote $\vec x = (x_1, x_2, \ldots, x_d)$, $\vec \omega =
(\omega_1, \ldots, \omega_d)$.

The total scheme remains much the same as in one dimension:

\begin{algorithm} \label{alg:twodim}{\sc Total Scheme in $d$ dimensions} \\
Input: signal $S$, the number of nonzero modes $B$ or its upper bound, accuracy factor $\epsilon$, success probability $1-\delta$, an upper bound of the signal energy $M$, the standard deviation of the white Gaussian noise $\sigma$.
\begin{enumerate}
    \item Initialize the representation signal $R$ to 0, set the maximum number of iterations $T=B\log(N)\log(\delta)/\epsilon^{2}$,.
    \item Test whether $\Vert S-R \Vert^2 \leq \iota \|R\|^2$. If yes, return the representation signal $R$ and the whole algorithm ends; else go to step 3.
    \item Locate Fourier Modes $\vec \omega$ for the signal $S-R$ by the new isolation and group test procedures.
    \item Estimate Fourier Coefficients at $\vec \omega$: $\widehat{(S-R)}(\vec \omega)$.
    \item If the total number of iterations is less than $T$, go to 2; else return the representation $R$.
  \end{enumerate}
\end{algorithm}

The most important modification with respect to the one dimensional case lies in the procedure to carry out step 3 of Algorithm \ref{alg:twodim}. We adapt the technique for
frequency identification to fit the high dimensional case;
it is given by the following procedure.

\begin{algorithm} {\sc Locate the Fourier mode in $d$ dimensions}\label{locate2d}
Input: signal $S$, accuracy factor $\epsilon$, success probability $1-\delta$, an upper bound of the signal energy $M$.
\begin{enumerate}
    \item Random permutations in d dimension.
    \item Isolate in one (arbitrarily picked) dimension $i$ to get a new signal $F(t)=S*H_k(t)$.
    \item For each dimension $i'$, find the $i'$th component $\vec \omega^*_{i'}$ of the significant frequency by Group Testing for the signal $F$ in the $i'$th dimension.
    \item Finally, estimate the Fourier coefficients in the frequency $\vec \omega =
    (\omega^*_0, \ldots, \omega^*_{d-1})$.
          Keep the frequency d-tuple if its Fourier coefficient is large.
\end{enumerate}
\end{algorithm}

Note that the computational cost of the above algorithm is quadratic in the number of dimensions. The
permutation involves a $d \times d$ matrix\footnote{ Note that generalizing to $d$ dimensions our 1-dimensional
practice of checking not only the central frequency found, but also nearby neighbors, would make this algorithm
exponential in $d$, which is acceptable for small $d$. For large $d$, we expect it would suffice to check 
a fixed number of randomly picked nearby neighbors, removing the exponential nature of this technical feature.}
The group test procedure in each dimension processes the {\it
same} isolation signal. If a filter with $B$ taps is used for the isolation, then it captures at least one
significant frequency in the pass region with probability $1/B$. The basic idea behind this procedure is that,
because of the sparseness of the Fourier representation, cutting the frequency domain into slices of width $1/B$
in 1 dimension, leaving the other dimensions untouched, leads to, with positive probability, a separation of the
important modes into different slices. After this essentially 1-dimensional isolation, we only need to identify
the coordinates of the isolated frequency mode. After isolation, we assume $F(\vec x)=A e^{2 \pi i \vec \omega
\cdot \vec x/N}$, where $A$ and $\vec \omega$ are unknown. To find $\omega_{j'}$, we sample in the $j'$-th
coordinate only, keeping $x_1, \ldots, x_{j'-1}, x_{j'+1},\ldots, x_d$ fixed, so that (for this step) $F(\vec
x)$ can be viewed as $A e^{2 \pi i \vec \omega \cdot \vec x/N}= \tilde{A}e^{2 \pi i \omega_{j'} x_{j'}/N} $,
where $\tilde{A}=Ae^{2 \pi i (x_1 \omega_1 + \ldots + x_{j'-1}\omega{j'-1}+x_{j'+1}\omega{j'+1}+ \ldots+x_d
\omega_d)}$, remains the same for different $x_{j'}$ and has the same absolute value as $A$, which we can do in
each dimension separately by the following argument.

If we just repeated the 1-dimensional technique in each dimension, that is, carried out isolation in each of the $d$ dimensions sequentially,
the time cost would be exponential in the dimension $d$. We discuss now in some detail the steps 1, 2, 3 of Algorithm \ref{locate2d}.

\subsubsection{Random Permutations}
\label{sect:permute}

In one dimensional RA$\ell$SFA, the isolation part includes random permutations and the construction of signals with one frequency dominant.
However, the situation is more complicated in higher dimensions, which is why we separated out the permutation step in the algorithm.

Recall that in one dimension, the signal is dilated and modulated randomly in order to
separate possibly neighboring frequencies. In higher dimensions, different modes can have identical coordinates in some of
the dimensions; they would continue to coincide in these dimensions if we just applied  ``diagonal'' dilation, i.e., if we carried
out dilation and modulation sequentially in the different dimensions. To separate such modes, we need to use
random matrices. We transform any point
$(x_1,x_2,\ldots, x_d)$ into $(y_1,\ldots, y_d)$ given by
\begin{equation}
 \left(\begin{array}{c}
y_1\\
\vdots \\
y_d \end{array}\right)=\left(\begin{array}{cccc}
a_{11} & a_{12}& \ldots & a_{1d}\\
\vdots & \vdots & \vdots& \vdots \\
a_{d1} & a_{d2}& \ldots & a_{dd}\end{array}\right)\left(\begin{array}{c}
x_1\\
\vdots \\
x_d \end{array}\right)+\left(\begin{array}{c}
b_1\\
\vdots \\
b_d \end{array}\right)
\end{equation}
where $A=\left (a_{ij}\right )$ is a random and invertible matrix, the $a_{ij}$ and the $b_i$ are chosen randomly, uniformly and independently,
and the arithmetic is modulo $N$.
For example, if $d=2, N=7, a_{11}=1, a_{12}=3, a_{21}=5, a_{22}=2, b_1=0, b_2=5 $, that is, 
\begin{equation}
 \left(\begin{array}{c}
y_1\\
y_2 \end{array}\right)=\left(\begin{array}{cc}
1 & 3\\
5 & 2 \end{array}\right)\left(\begin{array}{c}
x_1\\
x_2 \end{array}\right)+\left(\begin{array}{c}
0\\
5 \end{array}\right)
\end{equation}
the point $(1,2)$ gets mapped to
$(0,0)$, $(1,3)$ to $(3,2)$, and $(0,3)$ to $(2,4)$: even though points $(1,2)$ and $(1,3)$ have the same first coordinate, their
images don't share a coordinate; the same happens with points $(1,3)$ and $(0,3)$. For each dimension $i'$, the $i'$th components of frequencies are mapped by pairwise
independent permutations. Even adjacent points that differ in only one coordinate are destined to be separate with high
probability after these random permutations.

\subsubsection{Isolation}
\label{sect:twodimisolation}

After the random permutations, the high dimensional version of isolation can construct a sequence
$F_{0},F_{1},\ldots$ of signals, such that , for some j, $|\hat{F_{j}}(\omega^{'})|^{2} \geq 0.98\| F\|^{2}$.

\begin{algorithm} \label{alg:twodimiso}{\sc High Dimensional Isolation} \\
\text{}\hspace{10mm} Choose an arbitrary dimension $i$.
\begin{enumerate}
       \item    Filter on the dimension $i$ and leave all other dimensions alone, get the signal
      \begin{equation}
          F =S \star H_{k},
      \end{equation}
             where $H_{k}=\frac{\sqrt{N}}{2k+1}\chi_{[-k,k]}$ filters on the dimension $i$; the other dimensions are not affected.
      \item  Output new signals $F$ to be used in the Group Testing.
  \end{enumerate}
\end{algorithm}

\subsubsection{Group Testing for Each Dimension}
\label{sect:twodimgroup}

After the random permutation and isolation, we expect a $d$-dimensional signal with most of its energy
concentrated on one mode. The isolation
step effectively separates the $d$-dimensional frequency domain in a number of $d$-dimensional slices. Group
testing has to subdivide these slices.

One naive method is to apply $d$ dimensional filters in group testing, concentrating on $d$-dimensional cubic subregions in group testing that cover the whole area.
However, this leads to more cost. If the number of taps of this filter in one dimension is $2k+1$, we
obtain $(2k+1)^{d}$ subregions. Estimating the energies of all subregions slows down the total running time.
Consequently we instead locate each component of the significant frequency label
 separately. That is, we only use a filter to focus on one dimension and leave other dimensions alone.
 The energy of $2k+1$ regions are computed in every dimension. Hence, we
need to estimate the norm of $d(2k+1)$ intervals in total. This makes Group Testing linear in the number of dimensions,
instead of exponential as in the naive method.

Here is the procedure in Group Test:
\begin{algorithm}  \label{alg:twodimgroup}{\sc High Dimensional Group Test} \\
\text{}\hspace{10mm} For $i'=1,\ldots, d$
\begin{enumerate}
    \item Construct signals $\tilde{G}_j^{(i')}=F(t) * ( e^{2 \pi i j t_{i'}/(2l+1)} H_{l})$, $j=1,\ldots, 2l+1$,
where $H_l$ filters on $i'$th dimension and leave all other dimensions alone;
    \item Estimate and compare the energy of each $\tilde{G}_j^{(i')}$, $j=1,\ldots, 2l+1$, use the similar procedure in one
    dimensional group testing procedure. Find the candidates $\omega^{*}_{i'}$.
\end{enumerate}
\end{algorithm}

The reader may wonder how sampling works out for this $d$-dimensional algorithm. In Algorithm \ref{alg:twodimgroup}, we will need to
sample $\tilde{G}_j^{(i')} $ (which is the convolution of the (permuted version of) signal $S$ with 2 filters) to estimate its
energy; because filtering is done only in the $i'$-th dimension, we shall sample
$\tilde{G}_j^{(i')}(x_1, \ldots, x_{i'-1}, x_{i'}, x_{i'+1}, \ldots, x_d)$ for different $x_{i'}$, but keeping the other $x_j$ fixed, where $j \neq i'$.
The signal $F$ itself comes from the Isolation step, in which we filter in direction $i$, for which $S$ needs to be sampled, in this dimension only.
Together, for each choices $i'$ in Algorithm \ref{alg:twodimiso} and \ref{alg:twodimgroup}, this implies we have $(2k+1)\times(2l+1)$ different samples
of (the permuted version of) $S$, in which all but the $i$th coordinates of the samples $\vec x$ are identical.

\section{Conclusion}

We provide both theoretical and experimental evidence to support the advantage of the implementation of RA$\ell$SFA proposed here over the original one sketched in \cite{GGIMS}. Moreover, we extend RA$\ell$SFA to high dimensional cases. For functions with few, dominant Fourier modes, RA$\ell$SFA outperforms the FFT as $N$ increases. We expect that RA$\ell$SFA will be useful as a substitute for the FFT in potential applications that require processing such sparse signals or computing $B$-term approximations. This paper is just the beginning of a series of our papers and researches, many of which are in preparation. For example, the strong dependence of running time on the number of modes $B$ will be further lessened, and thus the algorithm would work for more interesting signals \cite{GMS}. Also, the application of RA$\ell$SFA in multiscale problems will be discussed in \cite{ZDR}.

\section*{Acknowledgments}
For discussions that were a great help, we would like to thank Bjorn Engquist, Weinan E, Olof Runborg, and Josko Plazonic.

\end{document}